\documentclass[3p]{elsarticle}

\usepackage{amsmath}
\usepackage{amsthm}
\usepackage{amssymb}
\usepackage{bm}
\usepackage{accents}
\usepackage{tikz}

\newtheorem{thm}{Theorem}

\newdefinition{rmk}{Remark}
\newproof{pf}{Proof}

\newcommand{\const}{\mathop{\rm const}\nolimits}

\numberwithin{equation}{section}

\graphicspath{{figs/}}

\journal{arXiv} 

\begin{document}

\begin{frontmatter}

\title{Subdomain solution decomposition method for nonstationary problems\tnoteref{label1}}
\tnotetext[label1]{The paper was prepared within the framework of state task No. FSRG-2021-0015.}

\author{P.N. Vabishchevich\corref{cor}\fnref{lab1,lab2}}
\ead{vabishchevich@gmail.com}
\cortext[cor1]{Correspondibg author.}

\address[lab1]{Nuclear Safety Institute, Russian Academy of Sciences,
              52, B. Tulskaya, 115191 Moscow, Russia}

\address[lab2]{North-Eastern Federal University, 58, Belinskogo st, Yakutsk, 677000, Russia}

\begin{abstract}
The reduction of computational costs in the numerical solution of nonstationary problems is achieved through splitting schemes.
In this case, solving a set of less computationally complex problems provides the transition to a new level in time.
The traditional construction approach of splitting schemes is based on an additive representation of the problem operator(s) and uses explicit-implicit approximations for individual terms.
Recently (Y. Efendiev, P.N. Vabishchevich. Splitting methods for solution decomposition in nonstationary problems. \textit{Applied Mathematics and Computation}. \textbf{397}, 125785, 2021), a new class of methods of approximate solution of nonstationary problems has been introduced based on decomposition not of operators but of the solution itself.
This new approach with subdomain solution selection is used in this paper to construct domain decomposition schemes.
The boundary value problem for a second-order parabolic equation in a rectangle with a difference approximation in space is typical.
Two and three-level schemes for decomposition of the domain with and without overlapping subdomains are investigated.
Our numerical experiments complement the theoretical results.
\end{abstract}


\end{frontmatter}

\section{Introduction}\label{sec:1}

Efficient computational algorithms for multidimensional boundary value problems are based on different domain decomposition methods (DDM) variants.
The theoretical foundations of iterative domain decomposition methods are most fully developed for elliptic boundary value problems \cite{ToselliWidlund2005,QuarteroniValli1999}.
Different variants of DDM with and without overlapping subdomains and with various exchange conditions at boundaries of subdomains (interfaces) are considered.

For approximate solutions to nonstationary problems, we focus primarily on applying unconditionally stable difference schemes \cite{LeVeque2007,Samarskii1989}.
For the new time level problem, we can use those or other options for the domain decomposition method for stationary problems \cite{QuarteroniValli1999,mathew2008domain}.
Taking into account the specificity of nonstationary problems allows us to construct (see, for example, the implementation based on the Schwartz \cite{cai1991additive,cai1994multiplicative}) an optimal iterative domain decomposition method,
where the number of iterations is independent of the discretization steps in time and space.

When targeting computational systems with massive parallelism, success is achieved by decomposing not only in space but also in time.
In this case, time is treated as an additional space dimension.
Various parallel algorithms are constructed for time-dependent problems: the waveform relaxation methods \cite{gander2007optimized}, space-time multigrid methods \cite{weinzierl2012geometric} and parareal algorithms \cite{gander2007analysis}.
For approximation by space, for example, standard continuous finite-element approximations are used, and by time, the discontinuous Galerkin method \cite{Thomee2006}.

The specificity of nonstationary problems is almost wholly taken into account using the iterative-free domain decomposition method.
For the approximate solution of boundary value problems for the second-order parabolic equation, there is a rapid decrease in the error in the boundary condition with increasing distance from the boundary.
When using the Schwarz alternating method with overlapping subdomains, without loss of accuracy of approximate solutions, we can confine ourselves to one iteration \cite{Kuznetsov1988}.
Decomposition of the domain is connected with the additive representation of the problem operator with a formulation of particular problems in subdomains.
On this methodological basis, we can construct different variants of splitting schemes.
Such regionally-additive schemes are introduced in our work \cite{vabishchevich1989difference}.
The study of regionally-additive schemes is conducted using the theory of stability (correctness) of additive operator-difference schemes \cite{Marchuk1990,VabishchevichAdditive}.
The focus of the construction of decomposition schemes is on:
The method of decomposition for a computational domain
The choice of decomposition operators (exchange boundary conditions)
The splitting scheme (approximation in time)
The critical results of the construction and study of unconditionally stable regionally-additive schemes are reflected in the book \cite{SamarskiiMatusVabischevich2002}.
The theory and practice of unsupervised domain decomposition methods for time-dependent problems are being developed in various directions (see, e.g., \cite{vabishchevich2008domain,vabishchevich2011substructuring,vabishchevich2018two,hansen2017additive,eisenmann2022variational}). 

In classical splitting schemes \cite{Marchuk1990,VabishchevichAdditive} we use an additive representation of the problem operator.
In this case, solving evolutionary problems for individual simpler operator terms ensures a less computationally costly transition to a new level in time.
In many nonstationary problems, it is more convenient to construct computationally acceptable subtasks based on solution decomposition, when simpler problems are formulated for individual solution components.
In particular, for the approximate solution of multiphysics problems, we can group the solution components (monophysics problems) for the corresponding system of coupled equations.
Such solution splitting schemes are proposed and investigated in \cite{efendiev2021splitting} for the approximate solution of the Cauchy problem in a finite-dimensional Hilbert space for a first-order evolution equation.
In \cite{vabishchevich2021solution}, solution splitting schemes are constructed for second-order evolution equations.
The principal difference between this approach from standard splitting schemes consists precisely in the primary decomposition of the solution, after which separate subproblems are formulated.

The critical point of the study is related to the solution decomposition, with the selection of individual parts of the solution.
For this purpose, the corresponding restriction operators are introduced.
Decomposition of the solution is most straightforward when solving nonstationary problems on a direct sum of subspaces.
The second acceptable variant of decomposition of a solution we associate \cite{vabishchevich2021solution} with an additive representation of a unit operator based on restriction and extension operators.
Computationally convenient splitting solution schemes are constructed using explicit-implicit time approximations when the corresponding operator matrix's diagonal part or triangular splitting is extracted.
The stability of two- and three-level splitting schemes is investigated using general results of the theory of stability (correctness) of operator-difference schemes \cite{Samarskii1989,SamarskiiMatusVabischevich2002}.

This paper aims to construct domain decomposition schemes for nonstationary problems based on a general solution splitting scheme \cite{efendiev2021splitting}.
We consider a typical boundary value problem for a second-order parabolic equation in a rectangle.
We implement a difference approximation over the space on a rectangular grid.
The additive representation of the solution is achieved by decomposing the computational domain into subdomains.
We construct restriction and extension operators with overlapping subdomains and without overlapping subdomains, thus providing an additive representation of a unit operator \cite{vabishchevich2021solution}.
A consequence of this is splitting the grid self-adjoint elliptic operator of the problem into a sum of self-adjoint elliptic operators.
These individual operators have uniquely defined problems in subdomains and corresponding exchange conditions between subdomains.
A system of coupled nonstationary equations for the solution components is formulated.
Two- and tree-level difference schemes of the first and second-order time approximation are used to solve the Cauchy problem.
Stability and convergence conditions have been established for the diagonal part of the system of equations operator matrix to the new time level.
The possibilities of using explicit-implicit approximations for the triangular splitting of the operator matrix are also considered.
The theoretical results are illustrated by numerical experiments for the parabolic model problem.

The paper is organized as follows.
Section \ref{sec:2} presents a boundary value problem for a second-order parabolic equation and its difference approximation by space.
Decomposition of the solution is discussed in Section \ref{sec:3}.
The restriction and extension operators are constructed when decomposition of the grid domain with and without overlapping subdomains.
Section \ref{sec:4} is central to our work.
Here we construct and study two- and three-level splitting schemes in various variants.
The three-level scheme of the second order of approximation by a time when the diagonal part of the operator matrix of our system of equations is taken to a new level by time is highlighted as the most interesting for practical use.
Numerical experiments for a model parabolic problem are given in Section \ref{sec:5}.
The results of our study are summarized in Section \ref{sec:6}.

\section{Problem formulation}\label{sec:2}

We illustrate the key elements of the technique for constructing splitting solution schemes based on the domain decomposition method on a model two-dimensional boundary value problem for a second-order parabolic equation.
Without limiting the generality, we will assume that the computational domain $\Omega$ is a rectangle, and standard difference approximations on a uniform rectangular grid are used.
The transition to more general multidimensional problems in irregular regions using the finite volume method and the finite element method with a lumped mass procedure is editorial.

Let 
\[
\Omega= \{\bm{x}\mid \bm{x}=(x_1,x_2),\ 0<x_{d}<l_{d}
\ d =1,2\} . 
\]
In $\Omega$, we solve the boundary value problem for the parabolic equation
\begin{equation}\label{2.1}
  \frac{\partial u}{\partial t} -  
  \sum_{d=1}^2 \frac \partial {\partial x_d}\left ( k(\bm{x})\frac{\partial u}{\partial x_d}\right ) + c(\bm x) u =
  f(\bm{x},t), \quad \bm{x} \in \Omega ,
  \quad 0 < t \leq T ,  
\end{equation}
under standard assumptions $0 < k_1 \leq  k(\bm{x}) \leq k_2,$ $0 < c_1 \leq c(\bm x) \leq c_2$.
We consider the case of setting homogeneous boundary conditions of the second kind:
\begin{equation}\label{2.2}
  k(\bm{x}) \frac{\partial u}{\partial n} = 0, 
  \quad \bm{x} \in \partial \Omega, 
  \quad 0 < t \leq T ,  
\end{equation}
where $n$ is the normal to $\partial\Omega$. 
The initial condition is given in the form
\begin{equation}\label{2.3}
 u(\bm{x},0) = u^0(\bm{x}), \quad  \bm{x}\in \Omega . 
\end{equation}

In the computational domain $\Omega$ we introduce a uniform grid $\omega$ with steps $h_1$ and $h_2$:
\[
\begin{split}
\omega=&\big\{\bm{x}\mid \bm{x} = x_{ij}= \big((i-0.5)h_1,(j-0.5)h_2 \big),
\quad i = 1,2,\dots,N_1,\\
&j = 1,2,\dots,N_2,
\quad N_d h_d = l_d,\quad d = 1,2\big\}.
\end{split}
\]
In our case, all grid nodes $ \omega$ lie inside the region.
For grid functions $y, v$ defined on $\omega$, let us introduce the finite-dimensional Hilbert space $H$ with the scalar product and norm
\[
(y,v) = \sum_{\bm{x}\in \omega} y(\bm{x})v(\bm{x}) h_1 h_2 ,
\quad \| v\| = (v,v)^{1/2} .
\]
For the operator $D=D^* > 0$ we define a Hilbert space $H_D$ with scalar product and norm $(y,v)_D = (D y,v), \ \|v\|_D = (Dv,v)^{1/2}$.

After approximating by space for the approximate solution $v(\bm x, t), \ \bm x \in \omega$ of the boundary value problem (\ref{2.1})--(\ref{2.3}) we have the Cauchy problem for the differential-operator equation:
\begin{equation}\label{2.4}
  \frac{d v}{d t} + A v = \varphi(t),
  \quad 0 < t \leq T,  
\end{equation} 
\begin{equation}\label{2.5}
  v(0) = v^0,
\end{equation} 
where $v(t) = v(\cdot,t)$ and, for example, $\varphi(\bm{x},t) = f(\bm{x},t), \  \bm{x} \in \omega$.

Given (\ref{2.1}) for the grid operator $A$ we can use the representation
\begin{equation}\label{2.6}
A =A^{(1)} + A^{(2)} .
\end{equation}
Here, the operators $A^{(1)}, A^{(2)}$ are associated with them
to the corresponding differential operator in one direction.
For all but the boundary nodes, the grid operator $A^{(1)}$
can be taken in the form
\[
\begin{split}
A^{(1)} v = & -
\frac{1}{h_1^2} k(x_1+0.5h_1,h_2) \big(v(x_1+h_1,h_2) - v(\bm{x})\big) \\
& + \frac{1}{h_1^2} k(x_1-0.5h_1,h_2) \big(v(\bm{x}) - v(x_1-h_1,h_2)\big) + \frac{1}{2} c(\bm x) v(\bm x), \\
& \qquad \bm{x} \in \omega,
\quad x_1 \neq 0.5h_1,
\quad x_1 \neq l_1-0.5h_1.
\end{split}
\]
In the boundary nodes, the approximation is carried out taking into account the boundary condition (\ref{2.2}):
\[
\begin{split}
A^{(1)} v = & -
\frac{1}{h_1^2} k(x_1+0.5h_1,h_2) \big(y(x_1+h_1,h_2) - y(\bm{x})\big) + \frac{1}{2} c(\bm x) v(\bm x), \\
& \qquad \bm{x} \in \omega,
\quad x_1 = 0.5h_1,
\end{split}
\]
\[
\begin{split}
A^{(1)} v =
& \ \frac{1}{h_1^2} k(x_1-0.5h_1,h_2) \big(y(\bm{x}) - y(x_1-h_1,h_2)\big) + \frac{1}{2} c(\bm x) v(\bm x), \\
& \qquad \bm{x} \in \omega,
\quad x_1 = l_1-0.5h_1.
\end{split}
\]
The grid operator $A^{(2)}$ is constructed similarly.

By direct calculations  \cite{Samarskii1989}, we make sure that
\[
  A^{(1)} = \big(A^{(2)}\big)^* \geq \frac{1}{2} c_1 I, 
  \quad A^{(2)} = \big(A^{(2)}\big)^* \geq \frac{1}{2} c_1 I, 
\]
where $I$ is a unit operator in $H$.
Given (\ref{2.6}), the operator $A$ is self-adjoint and non-negative in $H$:
\begin{equation}\label{2.7}
 A = A^* \geq \delta I ,
\end{equation} 
where $\delta = c_1$.
In addition, this grid operator approximates the corresponding differential operator with second order on $\tau$  for sufficiently smooth coefficients and solutions of the problem (\ref{2.1})--(\ref{2.3}).

For solution of the problem (\ref{2.4}), (\ref{2.5}), (\ref{2.7}), we have an a priori estimate
\begin{equation}\label{2.8}
\|v(t) \|_A^2 \leq \|v^0\|_A^2 + \frac{1}{2} \int_{0}^{t} \|\varphi(\theta)\|^2 d \theta .
\end{equation}
This estimate provides the stability of the solution concerning the initial data and the right-hand side.

Two-level schemes with weights are usually \cite{Samarskii1989,SamarskiiMatusVabischevich2002} used when approximating first-order evolution equations by time.
We assume, for simplicity, that the time grid is uniform: $t^n = n \tau, \ n = 0, 1, \ldots, N$ where $\tau = T/N$ is time step.
Let $v^n = v(\bm x, t^n), \ \ \bm x \in \omega , \ n = 0, 1, \ldots, N$ and
\[
t^{n+\sigma} = \sigma t^{n+1} + (1-\sigma) t^n,
\quad \varphi^{n+\sigma} = \sigma \varphi^{n+1} + (1-\sigma) \varphi^n ,
\]
where $\sigma = \const, \ \sigma \in [0,1]$ is the weight parameter.
For the approximate solution $y^{n} \approx v(t^{n}) , \ n = 0,1, \ldots, N,$ of the problem (\ref{2.4}), (\ref{2.5}), we use a difference scheme
\begin{equation}\label{2.9}
\frac{y^{n+1} - y^{n}}{\tau } + A \big(\sigma y^{n+1} + (1-\sigma) y^n \big) = \varphi^{n+\sigma} ,
\quad \ n = 0,1, \ldots, N-1 ,
\end{equation}
\begin{equation}\label{2.10}
y^0 = v^0 .
\end{equation}

\begin{thm}\label{t-1}
The scheme with weight (\ref{2.9}), (\ref{2.10}) is unconditionally stable at $\sigma \geq 1/2$. 
For the solution, the estimate 
\begin{equation}\label{2.11}
 \|y^{n+1}\|_A^2 \leq \| v^0 \|_A^2 + \frac{1}{2} \sum_{k=0}^{n} \tau \| \varphi^{k+\sigma}\|^2,
 \quad \ n = 0,1, \ldots, N-1 ,
\end{equation} 
is holds.
\end{thm}

\begin{pf}
We write (\ref{2.9}) as
\begin{equation}\label{2.12}
 (I + \tau D) \frac{y^{n+1} - y^{n}}{\tau } + A \frac{y^{n+1} + y^{n}}{2} = \varphi^{n+\sigma} ,
 \quad \ n = 0,1, \ldots, N-1 ,  
\end{equation} 
with the operator
\[
 D = \Big (\sigma - \frac{1}{2} \Big ) A .
\] 
At $\sigma \geq 1/2$ we have $D = D^* \geq 0$.
By multiplying (\ref{2.12}) scalarly in $H$ by $2 (y^{n+1} - y^{n})$, we get
\[
 2 \tau \Big \| \frac{y^{n+1} - y^{n}}{\tau } \Big \|^2 + 2 \big (D(y^{n+1} - y^{n}), y^{n+1} - y^{n} \big ) 
 + \|y^{n+1}\|_A^2 = \|y^{n}\|_A^2 + 2 \tau \Big (\varphi^{n+\sigma} , \frac{y^{n+1} - y^{n}}{\tau } \Big ) .
\]  
Given 
\[
\begin{split}
 \big (D(y^{n+1} - y^{n}, y^{n+1} - y^{n} \big ) & \geq 0, \\
 \Big (\varphi^{n+\sigma} , \frac{y^{n+1} - y^{n}}{\tau } \Big ) & \leq \Big \| \frac{y^{n+1} - y^{n}}{\tau } \Big \|^2
 + \frac{1}{4} \| \varphi^{n+\sigma}\|^2,
\end{split}
\] 
this gives the inequality
\[
 \|y^{n+1}\|_A^2 \leq  \|y^{n}\|_A^2 +  \frac{1}{2} \tau \| \varphi^{n+\sigma}\|^2 .
\] 
From this follows the provable estimate (\ref{2.11}).
\end{pf}

The stability estimate (\ref{2.11}) is the discrete analogue of the estimate (\ref{2.8}).
The scheme (\ref{2.9}), (\ref{2.10}) has the second order of accuracy on $\tau$ for smooth solutions of (\ref{2.4}), (\ref{2.5}) (symmetric scheme, Crank-Nicholson scheme) and the first-order for other values of the weight parameter $\sigma$.

To find the solution on the new time level, we solve the problem
\[
(I + \sigma \tau A) y^{n+1} = \psi^{n},
\quad \psi^{n} = \big(I - (1-\sigma) \tau A \big) y^{n+1} + \tau \varphi^{n+\sigma} ,
\quad \ n = 0,1, \ldots, N-1 .
\]
For this purpose, applying different variants of iterative methods of domain decomposition (block iterative methods) using computers of parallel architecture is natural.
The main specificity of nonstationary problems manifests in the fact that we have an excellent initial approximation $y^{n}$ for the unknown solution $y^{n+1}$.
Therefore, we can limit ourselves to a small number of iterations.
The most fully prominent features of nonstationary problems become apparent in constructing iteration-free domain decomposition methods.

\section{Decomposition of the solution}\label{sec:3}

We consider numerical methods for the approximate solution of the Cauchy problem (\ref{2.4}), (\ref{2.5}) based on a solution decomposition.
To simplify the problem on a new time level, instead of the space $H$, we work with the family of finite-dimensional Hilbert spaces $H_\alpha, \ \alpha = 1,\dots, p$; these spaces have significantly lower dimensionality than $H$.
We construct an additive representation of the desired solution from the original space $H$ from the solutions in these spaces.

For each of these spaces, the linear restriction operators $R_\alpha$ from the space $H$ to $H_\alpha$ are defined 
and the extension (interpolation or prolongation)  operators $R^*_\alpha$ from $H_\alpha$ to $H$:
\[
 R_\alpha: H \mapsto H_\alpha,
 \quad  R^*_\alpha: H_\alpha \mapsto H, 
 \quad \alpha = 1,\dots,p . 
\] 
The choice of the restriction and extension operators is subject to the condition (see \cite{vabishchevich2021solution,vabishchevich2017vector})
\begin{equation}\label{3.1}
 \sum_{\alpha =1}^{p} R_\alpha^* R_\alpha  = I .
\end{equation} 
If condition (\ref{3.1}) is satisfied, we have a representation for the solution $v \in H$ using $v_\alpha \in H_\alpha, \ \alpha = 1,\dots,p$:
\begin{equation}\label{3.2}
 v = \sum_{\alpha =1}^{p} R_\alpha^* R_\alpha v = \sum_{\alpha=1}^{p}  R_\alpha^* v_\alpha,
 \quad   v_\alpha = R_\alpha v,
 \quad \alpha = 1,\dots,p .  
\end{equation} 

Let us formulate the problem for finding $v_\alpha \in H_\alpha, \ \alpha = 1,\dots,p$ in the approximate solution of the Cauchy problem (\ref{2.4}), (\ref{2.5}).
Given (\ref{3.2}), equation (\ref{2.4}) is written as
\[
  \frac{d v}{d t} + A \sum_{\beta =1}^{p} R_\beta^* v_\beta = \varphi(t),
  \quad 0 < t \leq T .    
\] 
By multiplying this equation by $R_\alpha \ \alpha = 1,\dots,p,$ we get a system of equations
\begin{equation}\label{3.3}
  \frac{d v_\alpha}{d t} + \sum_{\beta =1}^{p} R_\alpha A R_\beta^* v_\beta = \varphi_\alpha(t),
 \quad \alpha = 1,\dots,p , 
  \quad 0 < t \leq T .  
\end{equation} 
For the right-hand sides, we have $ \varphi_\alpha(t) = R_\alpha \varphi(t), \ \alpha = 1,\dots,p$.  
The initial condition (\ref{2.5}) gives
\begin{equation}\label{3.4}
 v_\alpha(0) = v_\alpha^0, 
 \quad v_\alpha^0 = R_\alpha v^0 ,
 \quad \alpha = 1,\dots,p .  
\end{equation} 

It is convenient for us to write the system of equations (\ref{3.3}) in the form of one first-order equation for vector quantities.
Define a vector $\bm v = \{v_1, \ldots, v_p \} $ and $\bm \varphi = \{\varphi_1, \varphi_2, \ldots, \varphi_p \} $, and then (\ref{3.3}), (\ref{3.4}) gives us the Cauchy problem
\begin{equation}\label{3.5}
 \frac{d \bm v}{d t} + \bm A \bm v = \bm \varphi ,
\end{equation} 
\begin{equation}\label{3.6}
 \bm v(t) = \bm v^0 .
\end{equation} 
The operator matrix $\bm A$ is
\[
 \bm A = \{R_\alpha A R^*_\beta \} ,
 \quad \alpha, \beta  = 1,\dots,p . 
\] 

We consider the problem (\ref{3.5}), (\ref{3.6}) on the direct sum of spaces $\bm H = H_1  \oplus \ldots \oplus H_p$.
For $\bm v, \bm y \in \bm H$, the scalar product and norm are defined by the expressions
\[
 (\bm v, \bm y) = \sum_{1=1}^{p} (v_\alpha, y_\alpha)_\alpha ,
 \quad \|\bm v\| = (\bm v, \bm v)^{1/2} ,
\]
where $(v_\alpha , y_\alpha)_\alpha$ is the scalar product for $v_\alpha, y_\alpha  \in H_\alpha$ in $H_\alpha , \ \alpha  = 1,\dots, p$.  

We give a simple a priori estimate for the solution of (\ref{3.5}), (\ref{3.6}).
Multiplying equation (\ref{3.5}) scalarly in $\bm H$ by $d\bm v / dt$, we get
\begin{equation}\label{3.7}
 \Big \|\frac{d \bm v}{d t} \Big \|^2 + \frac{1}{2} \frac{d }{d t} (\bm A \bm v, \bm v) = \Big (\bm \varphi , \frac{d \bm v}{d t} \Big ) .
\end{equation} 
Given
\[
 \Big (\bm \varphi , \frac{d \bm v}{d t} \Big ) \leq \Big \|\frac{d \bm v}{d t} \Big \|^2 + \frac{1}{4} \|\bm \varphi\|^2 ,
\] 
from equality (\ref{3.7}) we get the estimate
\begin{equation}\label{3.8}
\big (\bm A \bm v(t), \bm v(t) \big) \leq  (\bm A \bm v^0, \bm v^0) + \frac{1}{2}\int_{0}^{t} \|\bm \varphi(\theta)\|^2 d \theta .
\end{equation} 
By considering (\ref{3.2}), we have  
\[
 (\bm A \bm v, \bm y) = \sum_{\alpha =1}^{p} \Big(\sum_{\beta =1}^{p} R_\alpha A R^*_\beta v_\beta , y_\alpha \Big)_\alpha  
 = \Big(A \sum_{\beta =1}^{p} R^*_\beta v_\beta , \sum_{\alpha =1}^{p} R^*_\alpha y_\alpha \Big) = (A v, y) .
\] 
Under the conditions of (\ref{2.7}), this gives
\begin{equation}\label{3.9}
 \bm A = \bm A^* > 0.
\end{equation} 
From the inequality (\ref{3.8}) follows an a priori estimate
\begin{equation}\label{3.10}
 \|\bm v(t)\|_{\bm A}^2 \leq  \|\bm v^0\|_{\bm A}^2 +  \frac{1}{2}\int_{0}^{t} \|\bm \varphi(\theta)\|^2 d \theta .
\end{equation} 

From the estimate (\ref{3.10}), we directly obtain the estimate (\ref{2.8}) for the solution of (\ref{2.4}), (\ref{2.5}), (\ref{2.7}).
It follows from
\[
 \|\bm \varphi\|^2 = \sum_{\alpha =1}^{p} (\varphi_\alpha, \varphi_\alpha)_\alpha = \sum_{\alpha =1}^{p} (R_\alpha \varphi, R_\alpha \varphi) 
 = \sum_{\alpha =1}^{p} ( R^*_\alpha R_\alpha \varphi, \varphi) = \|\varphi\|^2 ,
\] 
if the key property (\ref{3.1}) is taken into account.

The domain decomposition methods formulate a set of separate independent problems to find the solution components $v_\alpha , \ \alpha = 1,\dots, p,$ in the representation (\ref{3.2}).
The corresponding problems are defined by a particular choice of (i) restriction operators $R_\alpha, \ \alpha = 1,\dots, p,$ and (ii) time approximations for the approximate solution of (\ref{3.3}), (\ref{3.4}).

\begin{figure}[htbp]
\centering
\includegraphics[width=0.4\linewidth]{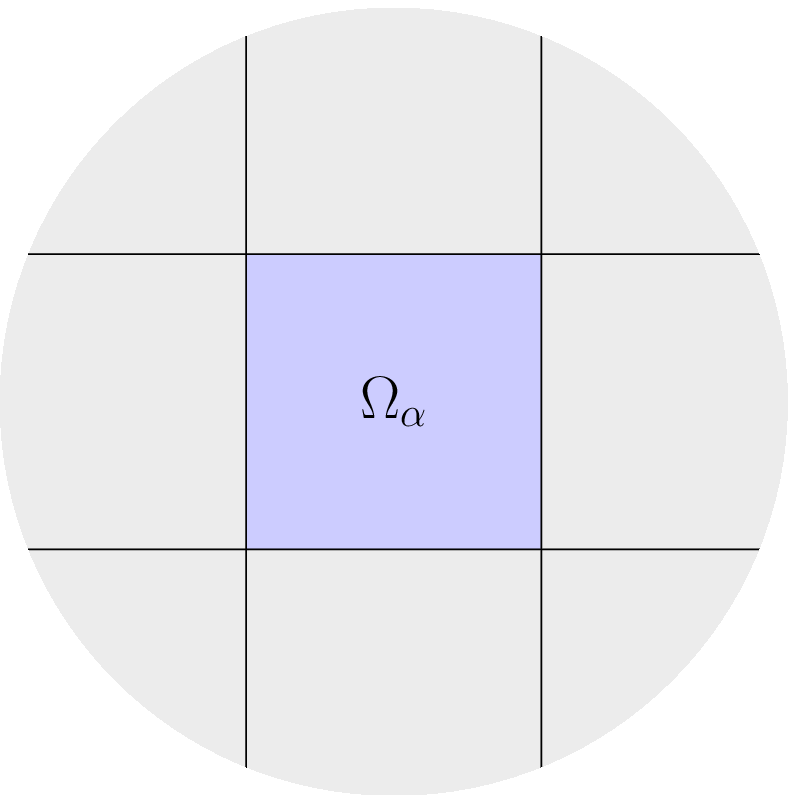} \hspace{0.1\linewidth} 
\includegraphics[width=0.4\linewidth]{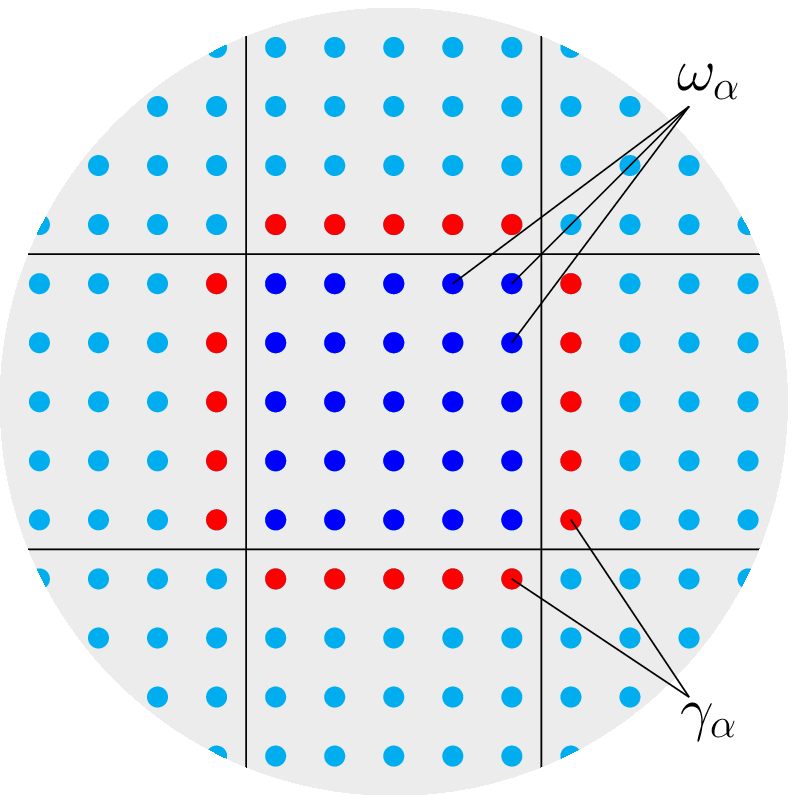} 
\caption{Decomposition of the domain (left) and grid (right) without overlapping subdomains.}
\label{f-1}
\end{figure}

We start by constructing the restriction operators in the decomposition of the domain
\[
\overline{\Omega} = \Omega \cup \partial \Omega = \bigcup_{\alpha =1}^{p} \overline{\Omega}_\alpha ,
\quad \overline{\Omega}_\alpha = \Omega_\alpha \cup \partial \Omega_\alpha ,
\quad \alpha = 1,\dots,p .
\]
In the multicoloring technique, each consolidated multi-connected subdomain $\overline{\Omega}_\alpha, \ \alpha = 1,\dots,p,$ consists of non-intersecting separate subdomains. This case ensures a parallel solution of the problem in the consolidated subdomain and minimizes the number of such subdomains $p$.
The simplest variant is characterized by non-overlapping subdomains: $\Omega_{\alpha \beta} = \Omega_\alpha \cap \Omega_\beta = \varnothing$ (see Fig.\ref{f-1}).
We assign to a separate subdomain $\Omega_\alpha$ part of the internal nodes $\omega_\alpha$ of the total computational grid $\omega$:
\[
\omega = \bigcup_{\alpha =1}^{p} \omega_\alpha ,
\quad \omega_\alpha = \{ \bm x \ | \ \bm x \in \omega, \ \bm x \in \Omega_\alpha \} ,
\quad \alpha = 1,\dots,p .
\]
The definition area of the grid operator $A$ for $Ay (\bm x), \ \bm x \in \omega_\alpha$ includes not only the inner grid nodes $\omega_\alpha$, but also the near-border nodes $\gamma_\alpha$ (Fig.\ref{f-1}). The near-border nodes lie in the neighboring subdomains of $\Omega_\alpha$; the data on these grid nodes are used to solve the problem in the $\Omega_\alpha$ subdomain.

When decomposing the domain without overlapping subdomains, the restriction operators for the grid function $v(\bm x), \ \bm x \in \omega$ are defined as follows:
\[
 v_\alpha (\bm x) = R_\alpha v (\bm x) = v(\bm x),
 \quad \bm x \in \omega_\alpha ,
 \quad \alpha = 1,\dots,p .    
\]  
The result of applying the extension operator $R^*_\alpha$ to the function $v_\alpha (\bm x), \ \bm x \in \omega_\alpha$ is the function
\[
 v (\bm x) = R^*_\alpha v_\alpha (\bm x) = \left \{ \begin{array}{ll}
  v_\alpha (\bm x),  & \bm x  \in \omega_\alpha,  \\
  0,  & \bm x  \in \omega \setminus \omega_\alpha,  
\end{array}
\right . 
 \quad \alpha = 1,\dots,p .    
\] 
Slightly more complex construction of the restriction and extension operators is used in the case of overlapping subdomains.

\begin{figure}[htbp]
\centering
\includegraphics[width=0.4\linewidth]{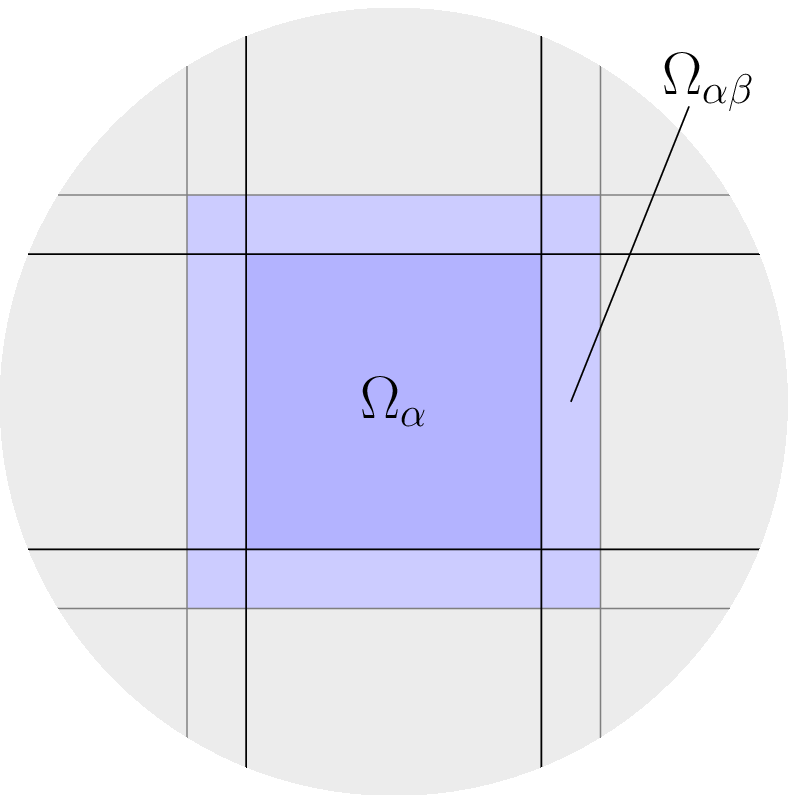} \hspace{0.1\linewidth} 
\includegraphics[width=0.4\linewidth]{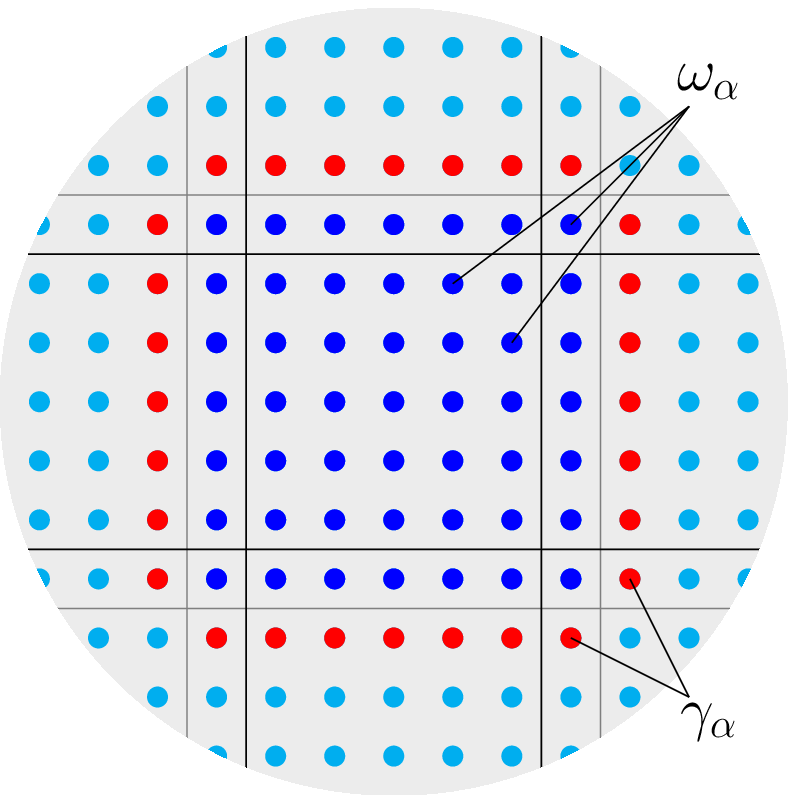} 
\caption{Decomposition of the domain (left) and grid (right) with overlapping subdomains.}
\label{f-2}
\end{figure}

The domain decomposition scheme with overlapping subdomains for a two-dimensional model problem on a uniform rectangular grid is shown in Fig.\ref{f-2}.
Compared to the case of decomposition without overlapping subdomains (compare with Fig.\ref{f-1}), the number of neighbors for a single subdomain and the number of common nodes increase.
We will construct the restriction and extension operators based on an auxiliary grid function $m(\bm x), \ \bm x \in \omega$.
For each node of the grid, we define the number of subdomains that this node belongs to:
\[
m(\bm x) = \sum_{\alpha =1}^{p} \chi_\alpha (\bm x),
\quad \chi_\alpha (\bm x) =
\left \{ \begin{array}{ll}
1, & \bm x \in \omega_\alpha, \\
0, & \bm x \in \omega \setminus \omega_\alpha,
\end{array}
\right .
\quad \bm x \in \omega .
\]
To decomposition in Fig.\ref{f-2} we have $m(\bm x) = 1,2,4, \ \ \bm x \in \omega$.
We set the restriction operators with the function $m(\bm x)$:
\[
v_\alpha (\bm x) = R_\alpha v (\bm x) = m^{-1/2}(\bm x) \, v(\bm x),
\quad \bm x \in \omega_\alpha ,
\quad \alpha = 1,\dots,p .
\]
For the extension operator $R^*_\alpha$, we have
\[
v (\bm x) = R^*_\alpha v_\alpha (\bm x) = m^{-1/2}(\bm x) \left \{ \begin{array}{ll}
v_\alpha (\bm x), & \bm x \in \omega_\alpha, \\
0, & \bm x \in \omega \setminus \omega_\alpha,
\end{array}
\right .
\quad \alpha = 1,\dots,p .
\]
The fulfillment of equality (\ref{3.1}) is checked directly.

\section{Splitting schemes}\label{sec:4} 

We construct time approximations for the Cauchy problem (\ref{3.5}), (\ref{3.6}), which provide a transition to a new time level by solving particular problems for the solution vector components.
First, we focus on splitting schemes \cite{VabishchevichAdditive}, which can be considered variants of standard two-level schemes with weights.

For an approximate solution of the problem (\ref{3.5}), (\ref{3.6}) we use the scheme
\begin{equation}\label{4.1}
 \frac{\bm y^{n+1} - \bm y^{n}}{\tau } + \bm A \big(\sigma \bm y^{n+1} + (1-\sigma) \bm y^n \big) = \bm \varphi^{n+\sigma} ,
 \quad   \ n = 0,1, \ldots, N-1 , 
\end{equation} 
\begin{equation}\label{4.2}
 \bm y^0 = \bm v^0 .
\end{equation} 

\begin{thm}\label{t-2}
The vector scheme with weight (\ref{4.1}), (\ref{4.2}) is unconditionally stable at $\sigma \geq 1/2$.
For the solution, the estimate
\begin{equation}\label{4.3}
\|\bm y^{n+1}\|_{\bm A}^2 \leq \|\bm v^0 \|_{\bm A}^2 + \frac{1}{2} \sum_{k=0}^{n} \tau \| \bm \varphi^{k+\sigma}\|^2,
\quad \ n = 0,1, \ldots, N-1 ,
\end{equation}
is holds.
For the approximate solution (\ref{2.4}), (\ref{2.5}), we have the estimate (\ref{2.11}).
\end{thm}

\begin{pf}
As in the proof of the theorem~\ref{t-1}, we write (\ref{4.1}) as
\begin{equation}\label{4.4}
 (\bm I + \tau \bm D) \frac{\bm y^{n+1} - \bm y^{n}}{\tau } + \bm A \frac{\bm y^{n+1} + \bm y^{n}}{2} = \bm \varphi^{n+\sigma} ,
 \quad \ n = 0,1, \ldots, N-1 ,  
\end{equation} 
where
\[
 \bm D = \Big (\sigma - \frac{1}{2} \Big ) \bm A .
\] 
Given (\ref{3.9}) and $\bm D = \bm D^* \geq 0$ at $\sigma \geq 1/2$, we get to the inequality
\[
 \|\bm y^{n+1}\|_{\bm A}^2 \leq \|\bm y^{n}\|_{\bm A}^2 + \frac{1}{2} \tau \|\bm \varphi^{n+\sigma}\|^2 .
\] 
From here, we have an estimate (\ref{4.3}).
Given the expressions for the solution and right-hand side components $v_\alpha^{n+1} , \varphi_\alpha^{n+\sigma}, \ \alpha = 1,2,\ldots, p$, from estimate (\ref{4.3}) follows (\ref{2.11}).  
\end{pf}

On the new time level, we solve the equation
\begin{equation}\label{4.5}
(\bm I + \sigma \tau \bm A) \bm y^{n+1} = \bm \psi^{n},
\quad \bm \psi^{n} = \big(\bm I - (1-\sigma) \tau \bm A \big) \bm y^{n+1} + \tau \bm \varphi^{n+\sigma} ,
\quad \ n = 0,1, \ldots, N-1 .
\end{equation}
For the individual solution components, we have a coupled system of equations
\[
y_\alpha^{n+1} + \sigma \tau \sum_{\beta =1}^{p} R_\alpha A R_\alpha^* y_\alpha^{n+1} = \psi_\alpha^{n} ,
\quad \alpha = 1,2, \ldots, p .
\]
To find the solution $y_\alpha^{n+1}, \ \alpha = 1,2, \ldots, p,$ we can use \cite{Saad2003} block Jacobi iterative methods and block triangular iterative methods (Seidel, SOR, SSOR).
Such methods are associated with additive and multiplicative Schwarz preconditioning in domain decomposition for approximate solutions of elliptic boundary value problems.
The grid problem (\ref{4.5}) is characterized by the presence of a small parameter --- the time step $\tau$.
In addition, there is a good initial approximation $y_\alpha^{n}, \ \alpha = 1,2, \ldots, p$.
These circumstances allow us to expect fast convergence of the corresponding iterative methods.
We consider iterative-free domain decomposition schemes for non-stationary problems.

To organize parallel computations, schemes for splitting the solution from the diagonal part of the operator matrix $\bm A$ on a new time level are of the most significant interest.
We will put
\[
\bm A_0 = \mathrm{diag} \{R_1 A R_1^*, R_2 A R_2^*, \ldots, R_p A R_p^* \} .
\]
For an approximate solution to the problem (\ref{3.5}), (\ref{3.6}) we use an explicit-implicit scheme of first-order accuracy
\begin{equation}\label{4.6}
 \frac{\bm w^{n+1} - \bm w^{n}}{\tau } + \bm A_0 \big(\sigma \bm w^{n+1} + (1-\sigma) \bm w^n \big) + (\bm A - \bm A_0) \bm w^{n}
 = \bm \varphi^{n+1} ,
 \quad   \ n = 0,1, \ldots, N-1 , 
\end{equation} 
\begin{equation}\label{4.7}
 \bm w^0 = \bm v^0 .
\end{equation} 

\begin{thm}\label{t-3}
At $\sigma \geq p/2$ the explicit-implicit scheme (\ref{4.6}), (\ref{4.7}) is unconditionally stable and the estimate
\begin{equation}\label{4.8}
 \|\bm w^{n+1}\|_{\bm A}^2 \leq \|\bm v^0 \|_{\bm A}^2 + \frac{1}{2} \sum_{k=0}^{n} \tau \| \bm \varphi^{k+1}\|^2,
 \quad \ n = 0,1, \ldots, N-1 ,
\end{equation} 
holds for the difference solution. 
\end{thm}

\begin{pf}
We write the scheme (\ref{4.6}) as
\begin{equation}\label{4.9}
 (\bm I + \tau \bm C) \frac{\bm w^{n+1} - \bm w^{n}}{\tau } + \bm A \frac{\bm w^{n+1} + \bm w^{n}}{2} = \bm \varphi^{n+1} ,
 \quad \ n = 0,1, \ldots, N-1 ,   
\end{equation} 
where
\[
 \bm C = \sigma \bm A_0 - \frac{1}{2} \bm A .
\] 
If $\bm C = \bm C^* \geq 0$, we have an estimate 
\[
 \|\bm w^{n+1}\|_{\bm A}^2 \leq \|\bm w^{n} \|_{\bm A}^2 + \frac{1}{2} \tau \| \bm \varphi^{n+1}\|^2,
 \quad \ n = 0,1, \ldots, N-1 , 
\] 
from which follows (\ref{4.8}).
We formulate conditions on the weight $\sigma$ to guarantee the non-negativity of the operator $\bm C$.

Given the inequality
\[
 \Big (\sum_{\alpha =1}^{p} a_\alpha \Big )^2 \leq p \sum_{\alpha =1}^{p} a_\alpha^2 ,
\] 
we have
\[
\begin{split}
 (\bm A_0 \bm v, \bm v) & = \sum_{\alpha =1}^{p} (R_\alpha A R_\alpha^* v_\alpha , v_\alpha )_\alpha 
 = \sum_{\alpha =1}^{p} \big ( (A^{1/2} R_\alpha^* v_\alpha )^2, 1 \big) \\
 & \geq \frac{1}{p} \Big ( \sum_{\alpha =1}^{p} A^{1/2} R_\alpha^* v_\alpha , \sum_{\beta =1}^{p} A^{1/2} R_\beta^* v_\beta \Big )
 = \frac{1}{p} (\bm A \bm v, \bm v) .
\end{split}
\] 
This results in the inequality
\begin{equation}\label{4.10}
 \bm A_0 \geq \frac{1}{p} \bm A . 
\end{equation} 
Under the constraints $\sigma \geq p/2$, the operator $\bm C$ is non-negative.
The theorem is proved.
\end{pf}

When applying domain decomposition schemes, the stability of the corresponding splitting schemes is essential.
The second key element of the study is related to the estimates of the error of the approximate solution, which are generated by using inhomogeneous approximations (explicit-implicit approximations) in time.
We briefly consider some vital elements of such a study.
We evaluate the accuracy of domain decomposition schemes by comparing the solution with the solution using the corresponding standard homogeneous time approximation schemes.

The error of the domain decomposition scheme (\ref{4.6}), (\ref{4.7}) is estimated by $\bm z^n = \bm w^n - \bm y^n, \ n = 1, 2, \ldots, N$.
Here $\bm y^n$ is the difference solution using the fully implicit ($\sigma =1$) scheme (\ref{4.1}), (\ref{4.2}).
From (\ref{4.1}), (\ref{4.2}) and (\ref{4.6}), (\ref{4.7}), we get
\begin{equation}\label{4.11}
\frac{\bm z^{n+1} - \bm z^{n}}{\tau } + \bm A \frac{\bm z^{n+1} + \bm z^{n}}{2 } = \bm \psi^{n+1} ,
\quad \ n = 0,1, \ldots, N-1 ,
\end{equation}
\begin{equation}\label{4.12}
\bm z^0 = 0 .
\end{equation}
The right-hand side (\ref{4.11}) is
\[
\bm \psi^{n+1} = \tau \bm D \frac{\bm y^{n+1} - \bm y^{n}}{\tau } - \tau \bm C \frac{\bm w^{n+1} - \bm w^{n}}{\tau } .
\]
For the difference scheme (\ref{4.11}), (\ref{4.12}), we have the estimate
\[
\|\bm z^{n+1}\|_{\bm A}^2 \leq \frac{1}{2} \sum_{k=0}^{n} \tau \| \bm \psi^{k+1}\|^2,
\quad \ n = 0,1, \ldots, N-1 .
\]
Taking into account the previously introduced notations for the right-hand side (\ref{4.11}), we obtain
\[
\begin{split}
\bm \psi^{n+1} & = \tau (\bm D - \bm C) \frac{d \bm v}{d t} (t^{n+1}) + \mathcal{O}(\tau^2)
= \tau (\bm A - \sigma \bm A_0) \frac{d \bm v}{d t} (t^{n+1}) + \mathcal{O}(\tau^2) \\
& = \bm \psi^{n+1}_0 + \bm \psi^{n+1}_1 + \mathcal{O}(\tau^2) .
\end{split}
\]
Part of the error
\[
\bm \psi^{n+1}_0 = \tau \bm A \frac{d \bm v}{d t} (t^{n+1}) = \tau A \frac{d v}{d t} (t^{n+1})
\]
is standard for schemes with weights.
For the error due to non-uniform time approximation, we have
\[
\bm \psi^{n+1}_1 = - \tau \sigma \bm A_0 \frac{d \bm v}{d t} (t^{n+1})
= - \tau \sigma \sum_{\alpha =1}^{p} R_\alpha^* A R_\alpha R_\alpha^* \frac{d v }{d t} (t^{n+1}) .
\]
Given $\|R_\alpha^*\| \leq 1, \ \ \alpha =1, 2, \ldots, p,$ this gives
\[
\|\bm \psi^{n+1}_1\| \leq \tau \sigma \sum_{\alpha =1}^{p} \Big \|A R_\alpha R_\alpha^* \frac{d v }{d t} (t^{n+1}) \Big \| .
\]
For our model problem (\ref{2.1}), (\ref{2.2}), when decomposing without overlapping subdomains, we obtain
\[
\sum_{\alpha =1}^{p} \Big \|A R_\alpha R_\alpha^* \frac{d v }{d t} (t^{n+1}) \Big \|
\leq p M \max_{\bm x \in \omega} \Big |\frac{d v }{d t} (t^{n+1}) \Big | ,
\quad M = \mathcal{O} (h^{-3/2}),
\quad h^2 = h_1^2 + h_2^2 .
\]
The error decreases as the overlap of subdomains are applied.
The main conclusion is that the domain decomposition scheme (\ref{4.6}), (\ref{4.7}) has an error $\mathcal{O}(\tau h^{-3/2})$.

In the second version of the domain decomposition schemes, we will relate to the triangular splitting of the operator matrix $\bm A$.
Let
\[
 \bm A_1 = \left (\begin{array}{cccc}
  \frac{1}{2} R_1 A R_1^* & 0 & \cdots & 0 \\
  R_2 A R_1^* & \frac{1}{2} R_2 A R_2^* & \cdots & 0 \\
  \cdots & \cdots & \cdots & 0 \\
   R_p A R_1^* & R_p A R_2^* & \cdots & \frac{1}{2} R_p A R_p^* \\
\end{array}
 \right ) ,
\quad
  \bm A_2 = \left (\begin{array}{cccc}
  \frac{1}{2} R_1 A R_1^* & R_1 A R_2^* & \cdots & R_1 A R_p^* \\
  0 & \frac{1}{2} R_2 A R_2^* & \cdots & R_2 A R_p^* \\
  0 & \cdots & \cdots & \cdots \\
  0 & 0 & \cdots & \frac{1}{2} R_p A R_p^* \\
\end{array}
 \right ) .
\]  
With this splitting, we have
\begin{equation}\label{4.13}
 \bm A = \bm A_1 + \bm A_2 ,
 \quad \bm A_1^* = \bm A_2 .
\end{equation} 

Among the schemes with triangular splitting (\ref{4.13}), we identify as the most exciting scheme of the alternating-triangular method \cite{Samarskii1989,VabishchevichAdditive} when solving problems with operators $\bm A_1$ and $\bm A_2$:
\[
\frac{\bm w^{n+1/2} - \bm w^{n}}{0.5 \tau } + \bm A_1 \bm w^{n+1/2} + \bm A_2 \bm w^n = \bm \varphi^{n+1/2} ,
\]
\[
\frac{\bm w^{n+1} - \bm w^{n+1/2}}{0.5 \tau } + \bm A_1 \bm w^{n+1/2} + \bm A_2 \bm w^{n+1/2} = \bm \varphi^{n+1/2} ,
\quad \ n = 0,1, \ldots, N-1 .
\]
It is convenient to write this operator analog of the Peaceman-Rachford scheme as a factorized scheme
\begin{equation}\label{4.14}
(\bm I + 0.5 \tau \bm A_1) (\bm I + 0.5 \tau \bm A_2) \frac{\bm w^{n+1} - \bm w^{n}}{\tau } + \bm A \bm w^n = \bm \varphi^{n+1/2} ,
\quad \ n = 0,1, \ldots, N-1 .
\end{equation}
The scheme (\ref{4.7}), (\ref{4.14}) has the second order of accuracy on $\tau$; this scheme can be considered as a one-iteration block SSOR method to implement the symmetric scheme ($\sigma = 1/2$ in (\ref{4.1}), (\ref{4.2})).

\begin{thm}\label{t-4}
The scheme of the alternating-triangular method  (\ref{4.7}), (\ref{4.13}), (\ref{4.14}) is unconditionally stable and the estimate 
\begin{equation}\label{4.15}
 \|\bm w^{n+1}\|_{\bm A}^2 \leq \|\bm v^0 \|_{\bm A}^2 + \frac{1}{2} \sum_{k=0}^{n} \tau \| \bm \varphi^{k+1/2}\|^2,
 \quad \ n = 0,1, \ldots, N-1 ,
\end{equation}
holds for the difference solution. 
\end{thm}

\begin{pf}
We write (\ref{4.14}) as
\[
(\bm I + \tau \bm C) \frac{\bm w^{n+1} - \bm w^{n}}{\tau } + \bm A \frac{\bm w^{n+1} + \bm w^{n}}{2} = \bm \varphi^{n+1/2} ,
\quad \ n = 0,1, \ldots, N-1 ,
\]
where
\[
\bm C = \frac{\tau }{4} \bm A_1 \bm A_2 .
\]
With splitting (\ref{4.13}) we have $\bm C = \bm C^* \geq 0$.
This property guarantees that the scheme (\ref{4.7}), (\ref{4.14}) is stable and that the estimate (\ref{4.15}) holds.
\end{pf}

The domain decomposition scheme (\ref{4.7}), (\ref{4.14}) is compared to the symmetric ($\sigma = 1/2$) scheme (\ref{4.1}), (\ref{4.2}).
For $\bm z^n = \bm w^n - \bm y^n, \ n = 1, 2, \ldots, N$ we have the problem (\ref{4.11}), (\ref{4.12}) with
\[
\bm \psi^{n+1/2} = - \frac{1 }{4} \tau^2 \bm A_1 \bm A_2 \frac{\bm w^{n+1} - \bm w^{n}}{\tau } .
\]
Thus, we get
\[
\bm \psi^{n+1/2} = \bm \psi^{n+1}_1 + \mathcal{O}(\tau^2) ,
\quad \bm \psi^{n+1/2}_1 = - \frac{1 }{4} \tau^2 \bm A_1 \bm A_2 \frac{d \bm v}{d t}(t^{n+1/2} ) .
\]
For $\bm \psi^{n+1/2}$ there is an estimate
\[
\| \bm \psi^{n+1/2}_1\| \leq \frac{p}{4} \tau^2 \sum_{\alpha, \beta =1}^{p} \Big \|A R_\alpha R_\alpha^* A R_\beta R_\beta^* \frac{d v }{d t} (t^{n+1/2}) \Big \| .
\]
For decomposing without overlapping subdomains, the central part of the error is estimated by the value
\[
\Big \|A R_\alpha R_\alpha^* A R_\alpha R_\alpha^* \frac{d v }{d t} (t^{n+1}) \Big \|
\leq M \max_{\bm x \in \omega} \Big |\frac{d v }{d t} (t^{n+1}) \Big | ,
\quad M = \mathcal{O} (h^{-7/2}) .
\]
Hence, the domain decomposition scheme (\ref{4.7}), (\ref{4.14}) has an error $\mathcal{O}(\tau^2 h^{-7/2})$.

The domain decomposition scheme (\ref{4.7}), (\ref{4.14}) based on the triangular division (\ref{4.13}) for the parabolic problem (\ref{2.1}), (\ref{2.2}) is inferior in accuracy to the domain decomposition scheme (\ref{4.6}), (\ref{4.7}) when the diagonal part $\bm A_0$ of the operator matrix $\bm A$ is brought to the top level in time. Moreover, the domain decomposition scheme (\ref{4.6}), (\ref{4.7}) has better prospects for parallelization similar to the block Jacobi method before the block Seidel method.

We can improve the accuracy of the domain decomposition scheme with the diagonal matrix $\bm A_0$ by using three-level time approximations.
The primary three-level scheme for the problem (\ref{3.5}), (\ref{3.6}) is the second-order accuracy scheme with a weight
\begin{equation}\label{4.16}
\frac{\bm y^{n+1} - \bm y^{n-1}}{2\tau } + \bm A \big(\sigma \bm y^{n+1} + (1-2\sigma) \bm y^n + \sigma \bm y^{n-1} \big) = \bm \varphi^{n} ,
\quad \ n = 1,2, \ldots, N-1 ,
\end{equation}
\begin{equation}\label{4.17}
\bm y^0 = \bm v^0 ,
\quad \bm y^1 = \widetilde{\bm v}^1 .
\end{equation}
To specify the second initial condition (\ref{4.17}) with the second-order accuracy, a two-level scheme is used:
\[
\frac{\widetilde{\bm v}^1 - \bm v^0}{\tau} + \bm A \frac{\widetilde{\bm v}^1 + \bm v^0}{2} = \bm \varphi^{1/2} .
\]

\begin{thm}\label{t-5}
The three-level scheme with weight (\ref{4.16}), (\ref{4.17}) is unconditionally stable at $\sigma \geq 1/4$. 
For the solution, the estimate 
\begin{equation}\label{4.18}
\begin{split}
 \Big (\sigma - \frac{1}{4} \Big ) \tau^2  & \Big \|\frac{\bm y^{n+1} - \bm y^{n}}{\tau } \Big \|_{\bm A}^2 
 + \Big \|\frac{\bm y^{n+1} + \bm y^{n}}{2 } \Big \|_{\bm A}^2 \\
 & \leq 
 \Big (\sigma - \frac{1}{4} \Big ) \tau^2  \Big \|\frac{\widetilde{\bm v}^1 - \bm v^0}{\tau } \Big \|_{\bm A}^2 
 + \Big \|\frac{\widetilde{\bm v}^1 + \bm v^0}{2 } \Big \|_{\bm A}^2
 + \frac{1}{2} \sum_{k=1}^{n} \tau \| \bm \varphi^{k}\|^2,
 \quad n = 0,1, \ldots, N-1 ,
\end{split}
\end{equation} 
is holds.
\end{thm}

\begin{pf}
It is convenient for us to introduce grid functions
\[
\bm s^n = \frac{1}{2} (\bm y^n + \bm y^{n-1}) ,
\quad \bm r^n = \frac{\bm y^n - \bm y^{n-1}}{\tau} .
\]
Given the identity
\[
\bm y^n = \frac{1}{4} (\bm y^{n+1} + 2 \bm y^n + \bm y^{n-1}) - \frac{1}{4} (\bm y^{n+1} - 2 \bm y^n + \bm y^{n-1}) ,
\]
let's write (\ref{4.16}) as
\begin{equation}\label{4.19}
\frac{\bm r^{n+1} + \bm r^{n}}{2 } + \bm D \frac{\bm r^{n+1} - \bm r^{n}}{\tau } + \bm A \frac{\bm s^{n+1} + \bm s^n}{2 } = \bm \varphi^{n} ,
\quad \ n = 1,2, \ldots, N-1 ,
\end{equation}
where
\[
\bm D = \Big (\sigma - \frac{1}{4} \Big ) \tau^2 \bm A .
\]
Multiplying equation (\ref{4.19}) scalarly in $\bm H$ by $2(\bm s^{n+1} - \bm s^{n}) = \tau (\bm r^{n+1} + \bm r^{n})$, we get
\[
\frac{\tau }{2} \|\bm r^{n+1} + \bm r^{n}\|^2 + (\bm D \bm r^{n+1}, \bm r^{n+1}) + \|\bm s^{n+1}\|_{\bm A}^2
\leq (\bm D \bm r^{n}, \bm r^{n}) + \|\bm s^{n}\|_{\bm A}^2 + \tau (\bm \varphi^{n} , \bm r^{n+1} + \bm r^{n}) .
\]
Since
\[
(\bm \varphi^{n} , \bm r^{n+1} + \bm r^{n}) \leq \frac{1}{2} \|\bm r^{n+1} + \bm r^{n}\|^2 + \frac{1}{2} \|\bm \varphi^{n} \|^2 ,
\]
we get
\[
(\bm D \bm r^{n+1}, \bm r^{n+1}) + \|\bm s^{n+1}\|_{\bm A}^2
\leq (\bm D \bm r^{n}, \bm r^{n}) + \|\bm s^{n}\|_{\bm A}^2 + \frac{\tau }{2} \|\bm \varphi^{n} \|^2 .
\]
At $\sigma \geq 1/4$, this gives
\[
(\bm D \bm r^{n+1}, \bm r^{n+1}) + \|\bm s^{n+1}\|_{\bm A}^2
\leq (\bm D \bm r^{1}, \bm r^{1}) + \|\bm s^{1}\|_{\bm A}^2 + \frac{1}{2} \sum_{k=1}^{n} \tau \| \bm \varphi^{k}\|^2 .
\]
Given the notations introduced, we have a provable estimate (\ref{4.18}).
\end{pf}

For $\sigma = 1/4$, we have an a priori estimate for the solution in semi-integer nodes when $\widetilde{\bm y}^{n+1/2} = 1/2 (\bm y^{n+1} + \bm y^{n})$.
In this case, for an approximate solution of the original problem (\ref{2.4}), (\ref{2.5}), from the estimate (\ref{4.18}), we obtain
\[
 \Big \|\frac{y^{n+1} + y^{n}}{2 } \Big \|_{A}^2 \leq 
 \Big \|\frac{\widetilde{v}^1 + v^0}{2 } \Big \|_{A}^2
 + \frac{1}{2} \sum_{k=1}^{n} \tau \| \varphi^{k}\|^2,
 \quad n = 0,1, \ldots, N-1 .
\] 

For the three-level scheme (\ref{4.16}), (\ref{4.17}), we compare an explicit-implicit three-level scheme with the allocation of the diagonal part of the operator matrix $\bm A$ on the new time level.
For the problem (\ref{3.5}), (\ref{3.6}), we use the scheme of second-order accuracy
\begin{equation}\label{4.20}
 \frac{\bm w^{n+1} - \bm w^{n-1}}{2\tau } + \bm A_0 \big(\sigma \bm w^{n+1} + (1-2\sigma) \bm w^n + \sigma \bm w^{n-1} \big) + (\bm A - \bm A_0) \bm w^{n}
 = \bm \varphi^{n} ,
 \quad   \ n = 1,2, \ldots, N-1 , 
\end{equation} 
\begin{equation}\label{4.21}
 \bm w^0 = \bm v^0 ,
 \quad \bm w^1 = \widetilde{\bm v}^1 .
\end{equation} 

\begin{thm}\label{t-6}
At $\sigma \geq p/4$ the three-level explicit-implicit scheme (\ref{4.20}), (\ref{4.21}) is unconditionally stable and for the solution there is the estimation
\begin{equation}\label{4.22}
\begin{split}
 \Big ( \bm C \frac{\bm y^{n+1} - \bm y^{n}}{\tau }, \frac{\bm y^{n+1} - \bm y^{n}}{\tau } \Big )
 & + \Big \|\frac{\bm y^{n+1} + \bm y^{n}}{2 } \Big \|_{\bm A}^2 \\
 \leq 
 \Big ( \bm C \frac{\widetilde{\bm v}^1 - \bm v^0}{\tau }, \frac{\widetilde{\bm v}^1 - \bm v^0}{\tau } \Big )
 & + \Big \|\frac{\widetilde{\bm v}^1 + \bm v^0}{2 } \Big \|_{\bm A}^2
 + \frac{1}{2} \sum_{k=1}^{n} \tau \| \bm \varphi^{k}\|^2, 
 \quad n = 0,1, \ldots, N-1 ,
\end{split}
\end{equation} 
where
\[
 \bm C = \sigma \tau^2 \bm A_0 - \frac{1}{4} \tau^2 \bm A \geq 0.
\]
\end{thm}

\begin{pf}
Let now
\[
\bm s^n = \frac{1}{2} (\bm w^n + \bm w^{n-1}) ,
\quad \bm r^n = \frac{\bm w^n - \bm w^{n-1}}{\tau} .
\]
Given this, we write (\ref{4.20}) as
\begin{equation}\label{4.23}
\frac{\bm r^{n+1} + \bm r^{n}}{2 } + \bm C \frac{\bm r^{n+1} - \bm r^{n}}{\tau } + \bm A \frac{\bm s^{n+1} + \bm s^n}{2 } = \bm \varphi^{n} ,
\quad \ n = 1,2, \ldots, N-1 .
\end{equation}
Considering the inequality (\ref{4.10}), we have $ \bm C \geq 0$.
From (\ref{4.23}), we obtain the inequality
\[
(\bm C \bm r^{n+1}, \bm r^{n+1}) + \|\bm s^{n+1}\|_{\bm A}^2
\leq (\bm C \bm r^{1}, \bm r^{1}) + \|\bm s^{1}\|_{\bm A}^2 + \frac{1}{2} \sum_{k=1}^{n} \tau \| \bm \varphi^{k}\|^2 .
\]
From this follows the estimate (\ref{4.22}).
\end{pf}

Let us compare the solution of the domain decomposition scheme (\ref{4.20}), (\ref{4.21}) and the solution of the scheme (\ref{4.16}), (\ref{4.17}) with $\sigma =1/4 \ (\bm D = 0)$.
We use the following notations
\[
\bm z^n = \bm w^n - \bm y^n,
\quad \bm s^n = \frac{1}{2} (\bm z^n + \bm z^{n-1}) ,
\quad \bm r^n = \frac{\bm z^n - \bm z^{n-1}}{\tau} ,
\quad n = 1, 2, \ldots, N .
\]
From (\ref{4.16}), (\ref{4.17}) and (\ref{4.20}), (\ref{4.21}), we have
\begin{equation}\label{4.24}
\frac{\bm r^{n+1} + \bm r^{n}}{2 } + \bm A \frac{\bm s^{n+1} + \bm s^{n}}{2 } = \bm \psi^{n} ,
\quad \ n = 0,1, \ldots, N-1 ,
\end{equation}
\begin{equation}\label{4.25}
\bm z^0 = 0 ,
\quad \bm z^1 = 0.
\end{equation}
For the right-hand side (\ref{4.24}), we use the representation
\[
\bm \psi^{n} = - \bm C \frac{\bm w^{n+1} - 2\bm w^{n} + \bm w^{n-1}}{\tau^2 } .
\]
Given homogeneous initial conditions for the difference scheme (\ref{4.24}), (\ref{4.25}), we have the estimate
\[
\Big \|\frac{\bm z^{n+1} + \bm z^{n}}{2} \Big \|_{\bm A}^2 \leq \frac{1}{2} \sum_{k=1}^{n} \tau \| \bm \psi^{k}\|^2,
\quad \ n = 1,2, \ldots, N-1 .
\]
For the right-hand side (\ref{4.24}), we obtain
\[
\begin{split}
\bm \psi^{n+1} & = - \bm C \frac{d^2 \bm v}{d t^2} (t^{n}) + \mathcal{O}(\tau^3)
= \tau^2 \Big (\frac{1}{4} \bm A - \sigma \bm A_0 \Big ) \frac{d^2 \bm v}{d t^2} (t^{n}) + \mathcal{O}(\tau^3) \\
& = \bm \psi^{n}_0 + \bm \psi^{n}_1 + \mathcal{O}(\tau^3) .
\end{split}
\]
For the standard part of the approximation error, we have
\[
\bm \psi^{n}_0 = \frac{1}{4} \tau^2 \bm A \frac{d^2\bm v}{d t^2} (t^{n}) = \frac{1}{4} \tau^2 A \frac{d^2 v}{d t^2} (t^{n}) .
\]
The error of the decomposition solution is
\[
\bm \psi^{n+1}_1 = - \sigma \tau^2 \bm A_0 \frac{d^2 \bm v}{d t^2} (t^{n})
= - \sigma \tau^2 \sum_{\alpha =1}^{p} R_\alpha^* A R_\alpha R_\alpha^* \frac{d v }{d t} (t^{n+1}) .
\]
and
\[
\|\bm \psi^{n+1}_1\| \leq \sigma \tau^2 \sum_{\alpha =1}^{p} \Big \|A R_\alpha R_\alpha^* \frac{d^2 v }{d t^2} (t^{n}) \Big \| .
\]
For the decomposition of the domain without overlapping subdomains, we have
\[
\sum_{\alpha =1}^{p} \Big \|A R_\alpha R_\alpha^* \frac{d^2 v }{d t^2} (t^{n}) \Big \|
\leq p M \max_{\bm x \in \omega} \Big |\frac{d^2 v }{d t^2} (t^{n}) \Big | ,
\quad M = \mathcal{O} (h^{-3/2}) .
\]
Thus the three-level domain decomposition scheme (\ref{4.20}), (\ref{4.21}) has an error $\mathcal{O}(\tau^2 h^{-3/2})$;
this scheme has significantly higher accuracy than the two-level domain decomposition scheme (\ref{4.6}), (\ref{4.7}).

\section{Numerical experiments}\label{sec:5}

We consider a test problem in unit square ($l_1 = l_2 = 1$) for the equation (\ref{2.1}) with
\[
k(\bm x) = 1,
\quad c(\bm x) = 1,
\quad f(\bm x,t) = \left \{ \begin{array}{ll}
25, & \bm x \in \Omega_f = [0.5,0.75] \times [0, 0.75] , \\
0, & \bm x \in \Omega \setminus \Omega_f .
\end{array}
\right .
\]
We find a solution with homogeneous initial conditions ($u^0 = 0$ in (\ref{2.3}) at $T=0.1$.
We compare the domain decomposition schemes with the standard schemes for the approximate solution of the problem (\ref{2.4}), (\ref{2.5}).
When using two-level schemes (see (\ref{2.9}), (\ref{2.10})), the primary one is the implicit scheme of first-order accuracy
\begin{equation}\label{5.1}
\frac{y^{n+1} - y^{n}}{\tau } + A y^{n+1} = \varphi^{n+1} ,
\quad \ n = 0,1, \ldots, N-1 ,
\end{equation}
\begin{equation}\label{5.2}
y^0 = 0 .
\end{equation}
In the class of three-level schemes (see (\ref{4.15}), (\ref{4.16})), we use a second-order accuracy scheme
\begin{equation}\label{5.3}
\frac{y^{n+1} - y^{n-1}}{2\tau } + \frac{1}{4} A (y^{n+1} + 2y^{n} + y^{n-1}) = \varphi^{n} ,
\quad \ n = 1,2, \ldots, N-1 .
\end{equation}
At the first time step, we have
\begin{equation}\label{5.4}
 \frac{y^{1} - y^{0}}{\tau } + \frac{1}{2} A (y^{1} + y^{0}) = \varphi^{1/2} . 
\end{equation} 

\begin{figure}[htbp]
\centering
\includegraphics[width=0.45\linewidth]{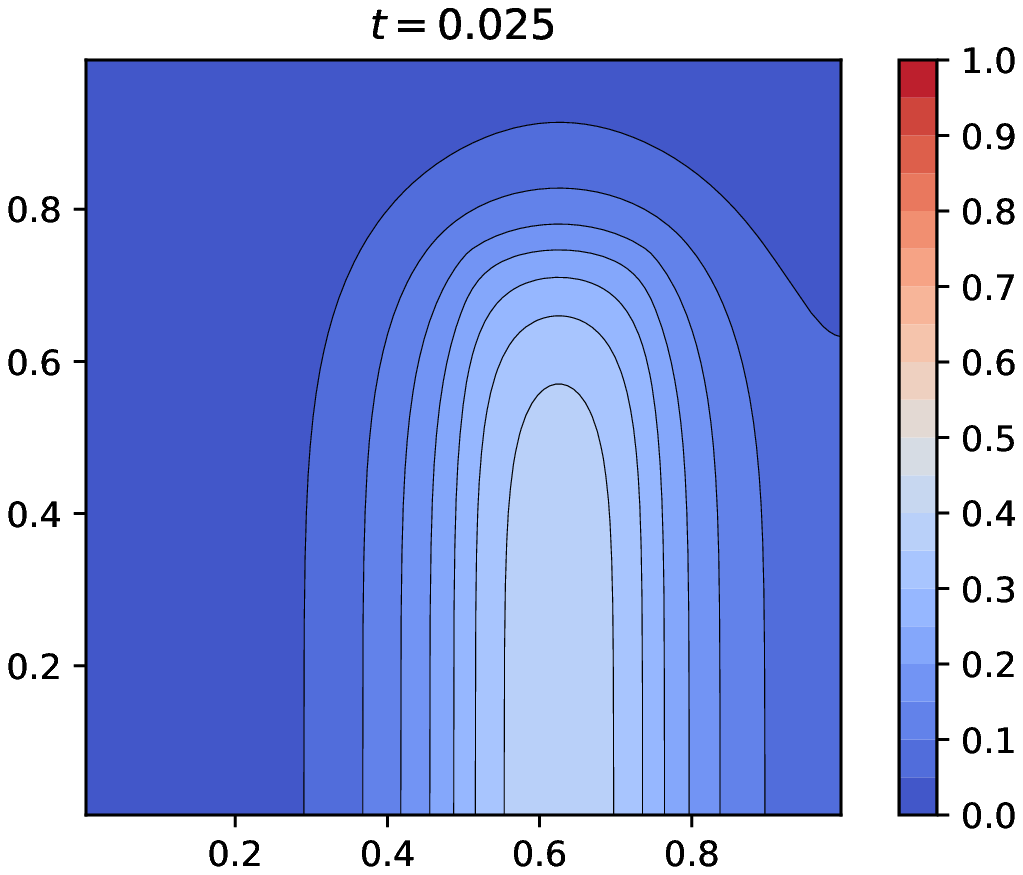} \includegraphics[width=0.45\linewidth]{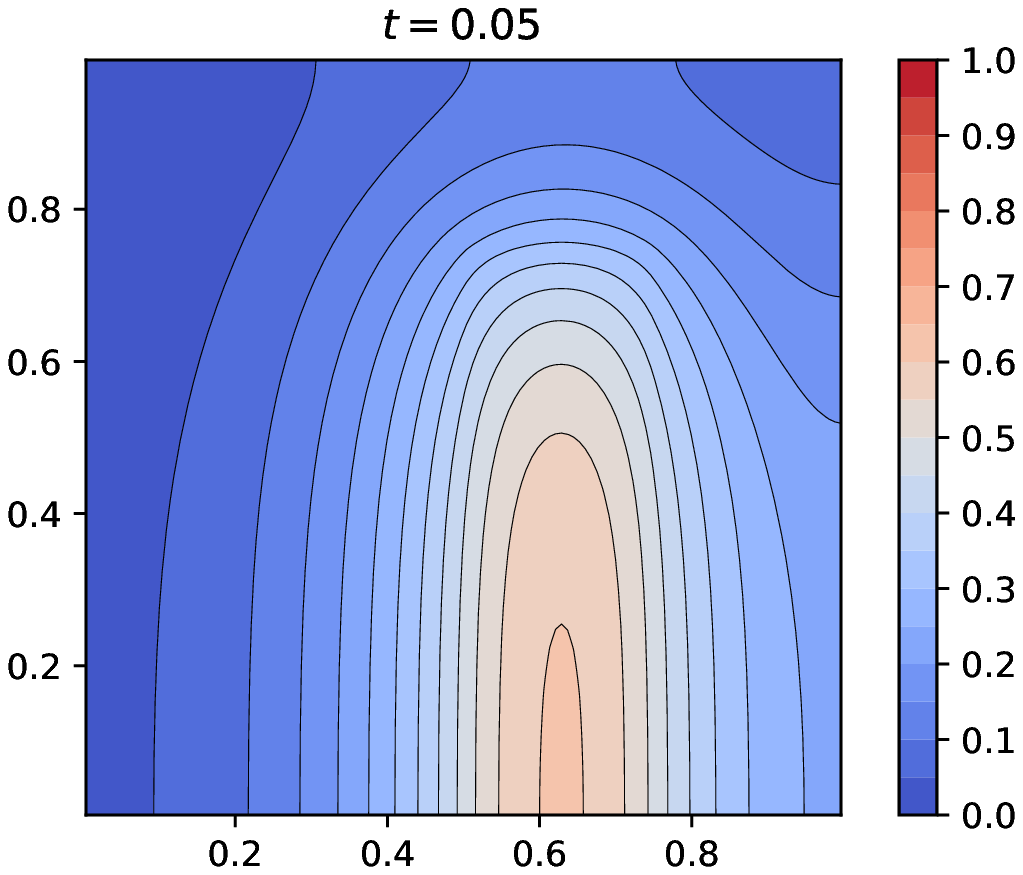} \\
\includegraphics[width=0.45\linewidth]{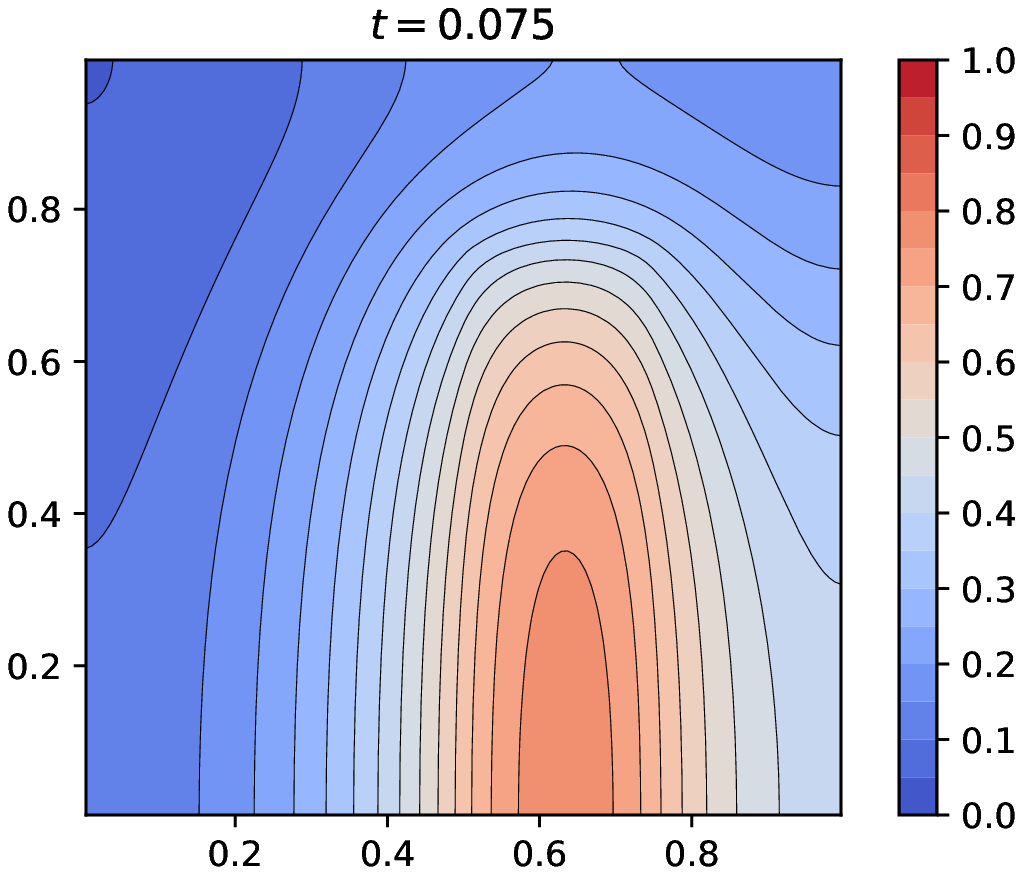} \includegraphics[width=0.45\linewidth]{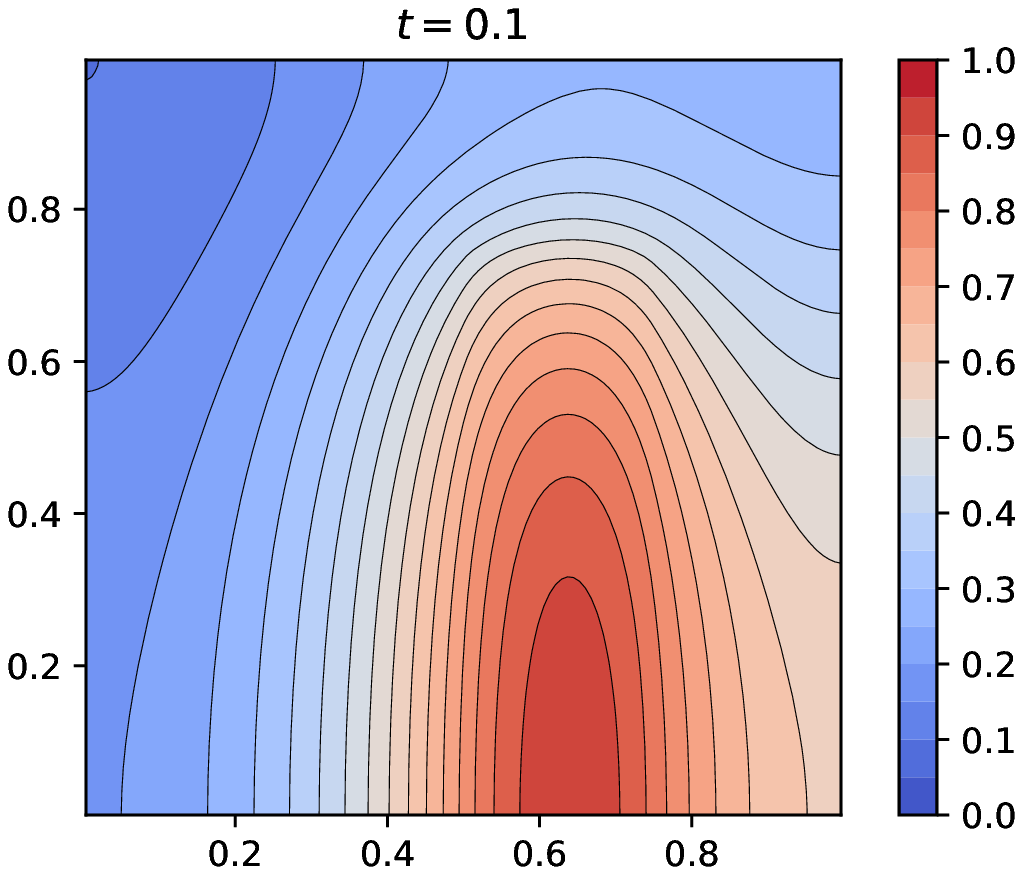} \\
\caption{The solution to specific points in time: $N_1 = N_2 = 128$.}
\label{f-3}
\end{figure}

\begin{figure}[htbp]
\centering
\includegraphics[width=0.45\linewidth]{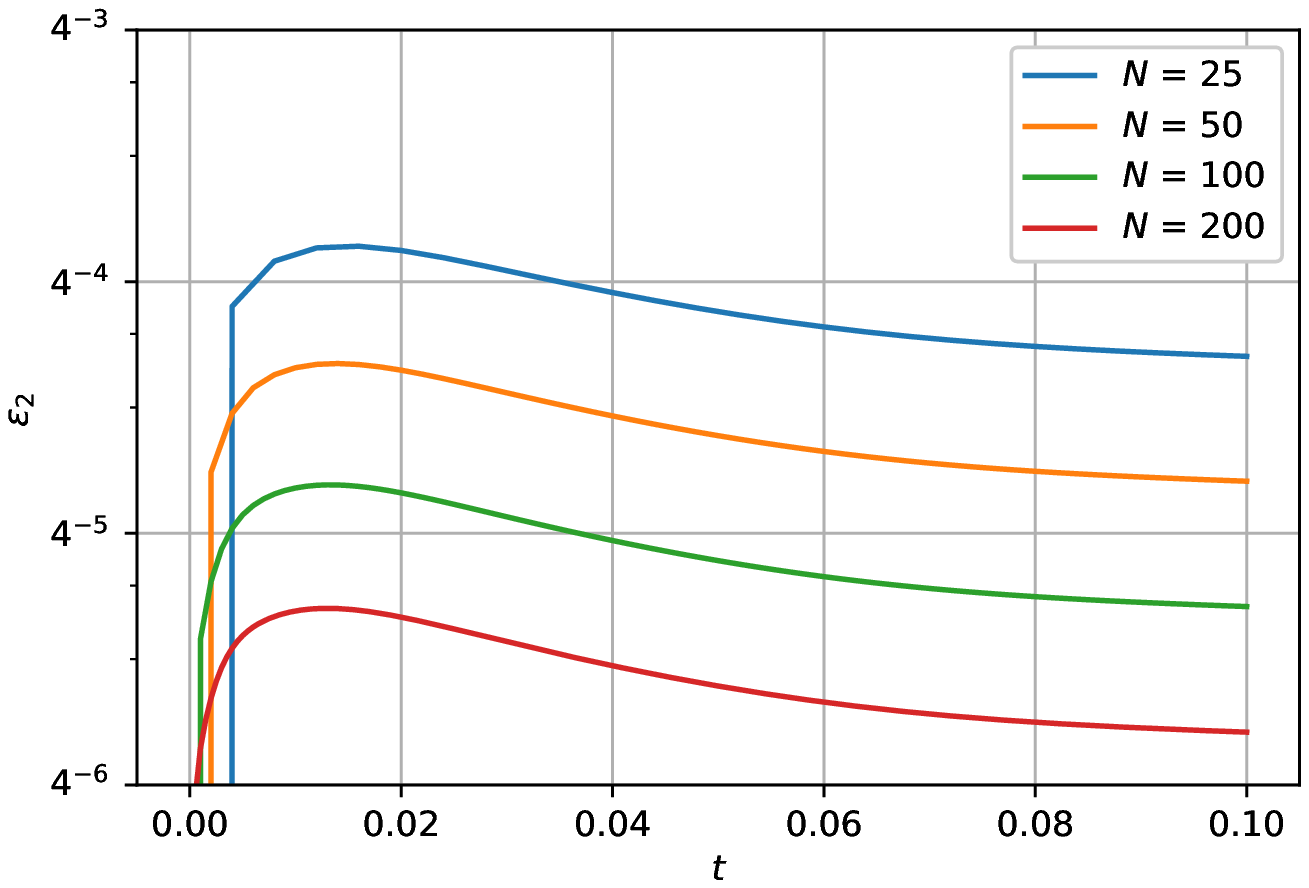} \includegraphics[width=0.45\linewidth]{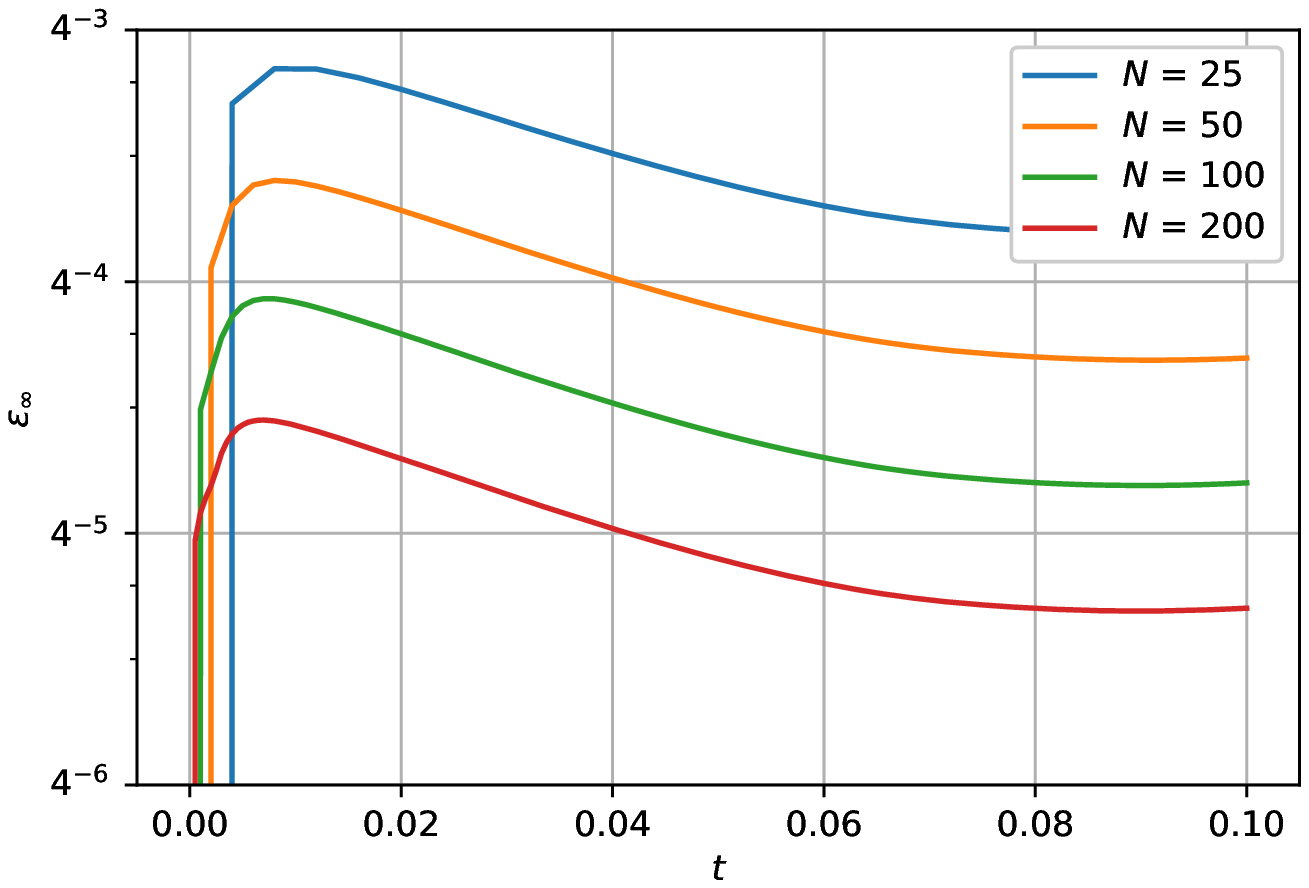} \\
\includegraphics[width=0.45\linewidth]{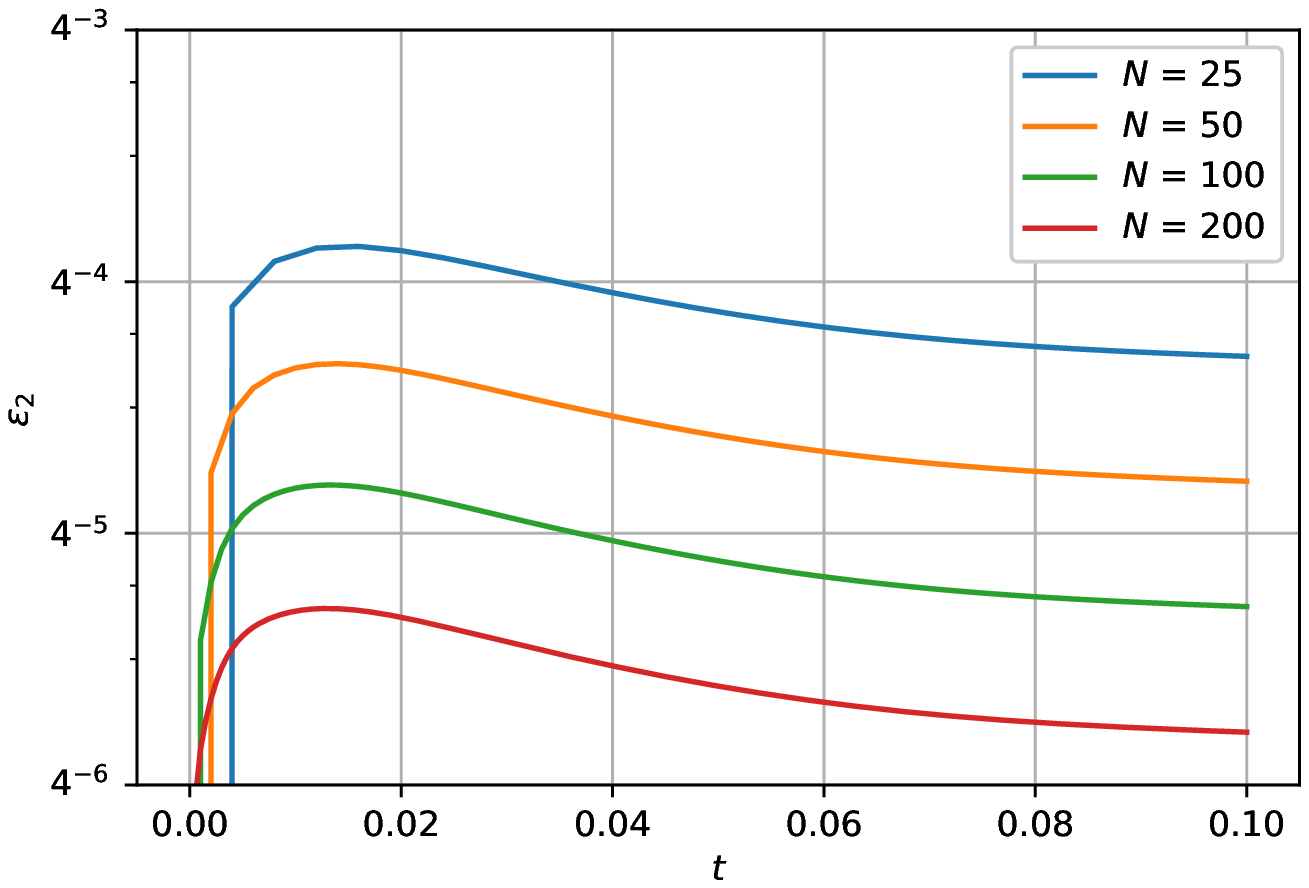} \includegraphics[width=0.45\linewidth]{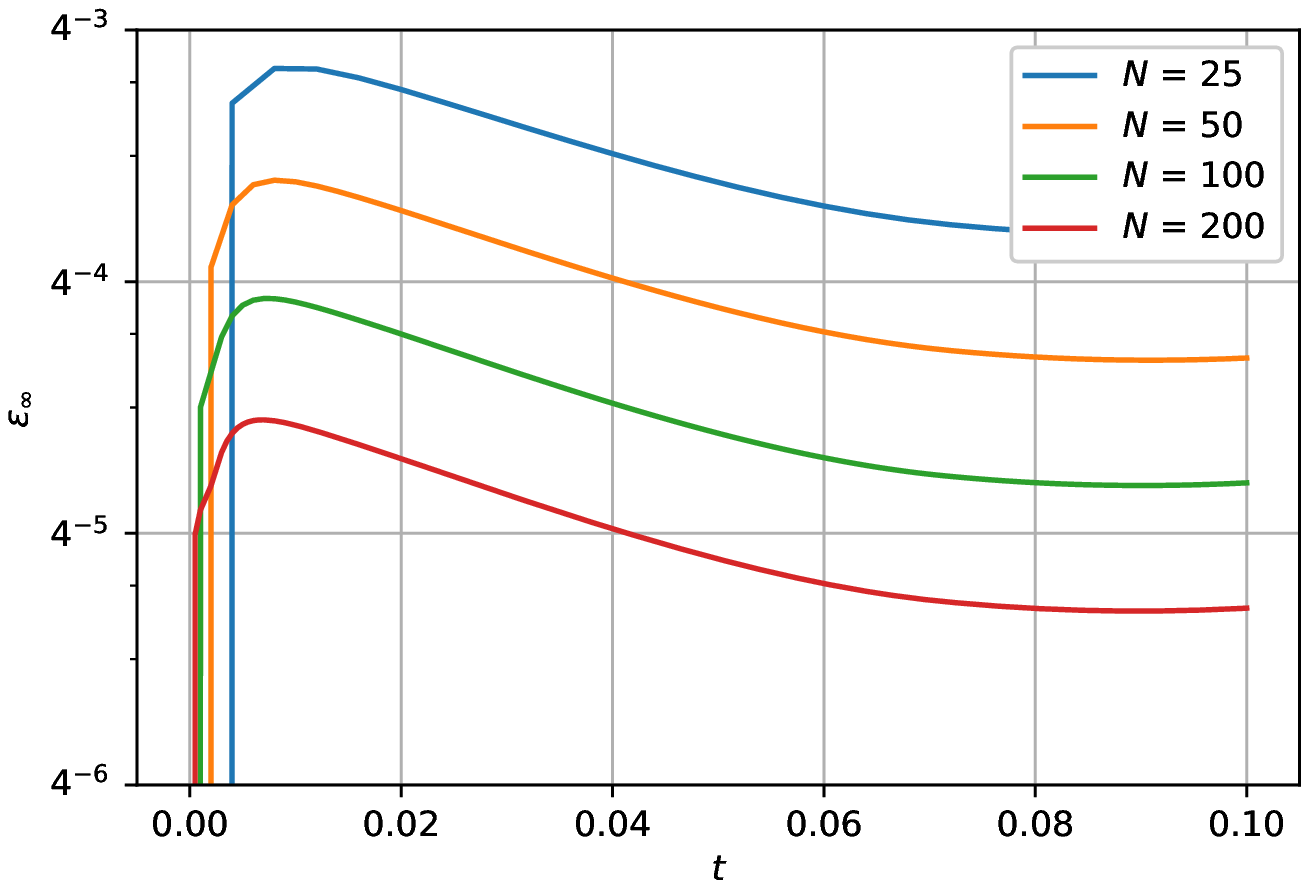} \\
\caption{The accuracy of the two-level scheme: $N_1 = N_2 = 64$ (top), $N_1 = N_2 = 128$ (bottom).}
\label{f-4}
\end{figure}

\begin{figure}[htbp]
\centering
\includegraphics[width=0.45\linewidth]{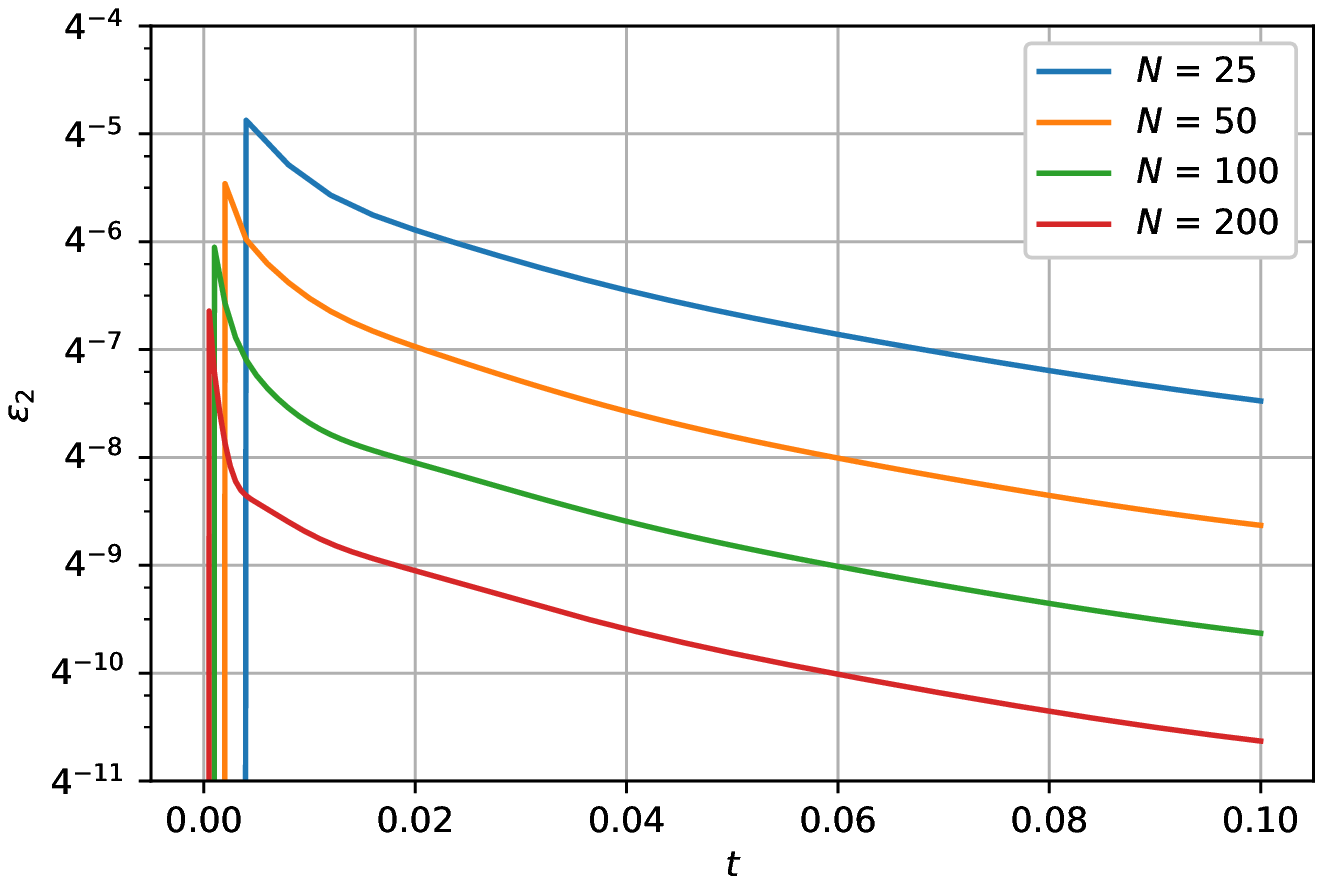} \includegraphics[width=0.45\linewidth]{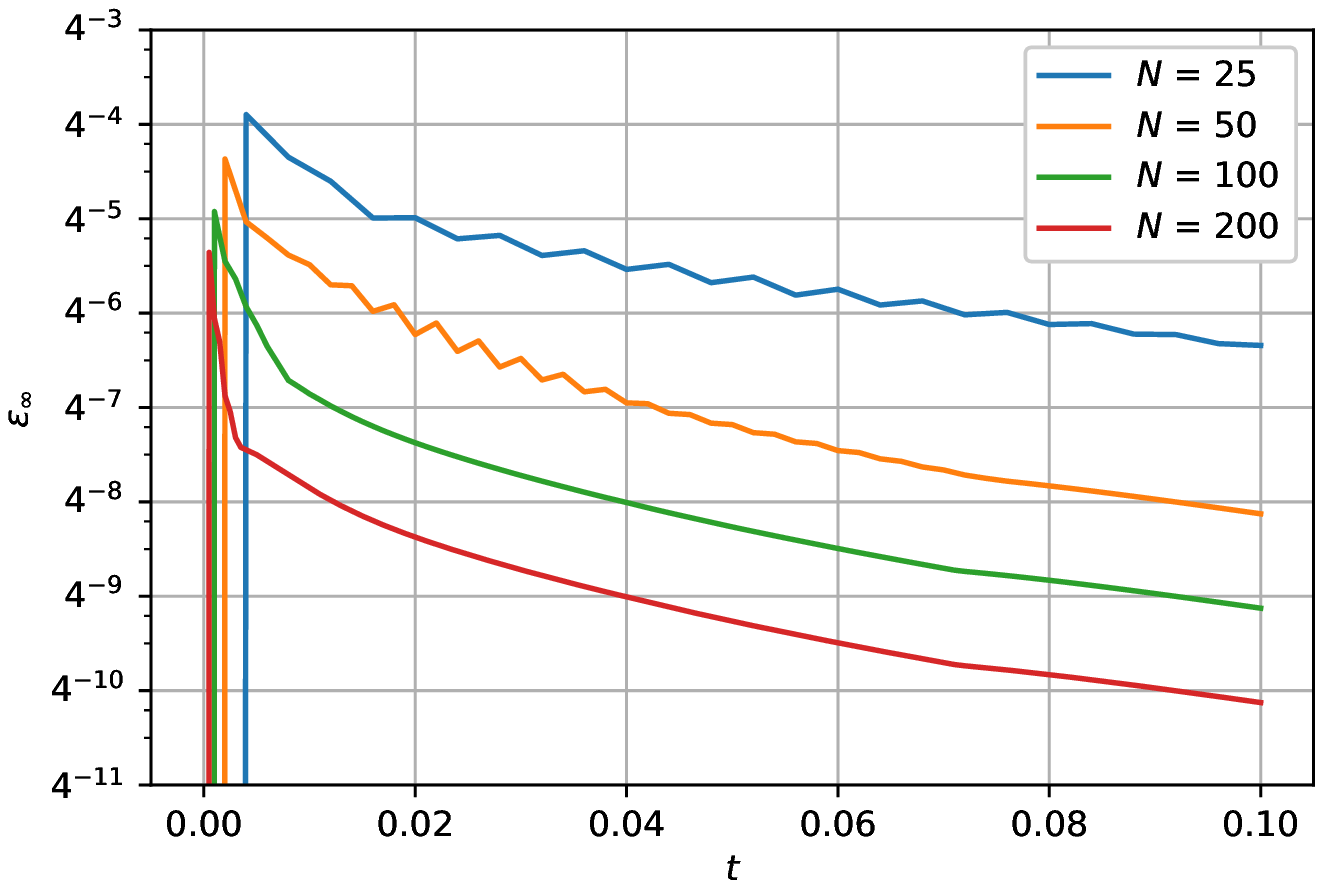} \\
\includegraphics[width=0.45\linewidth]{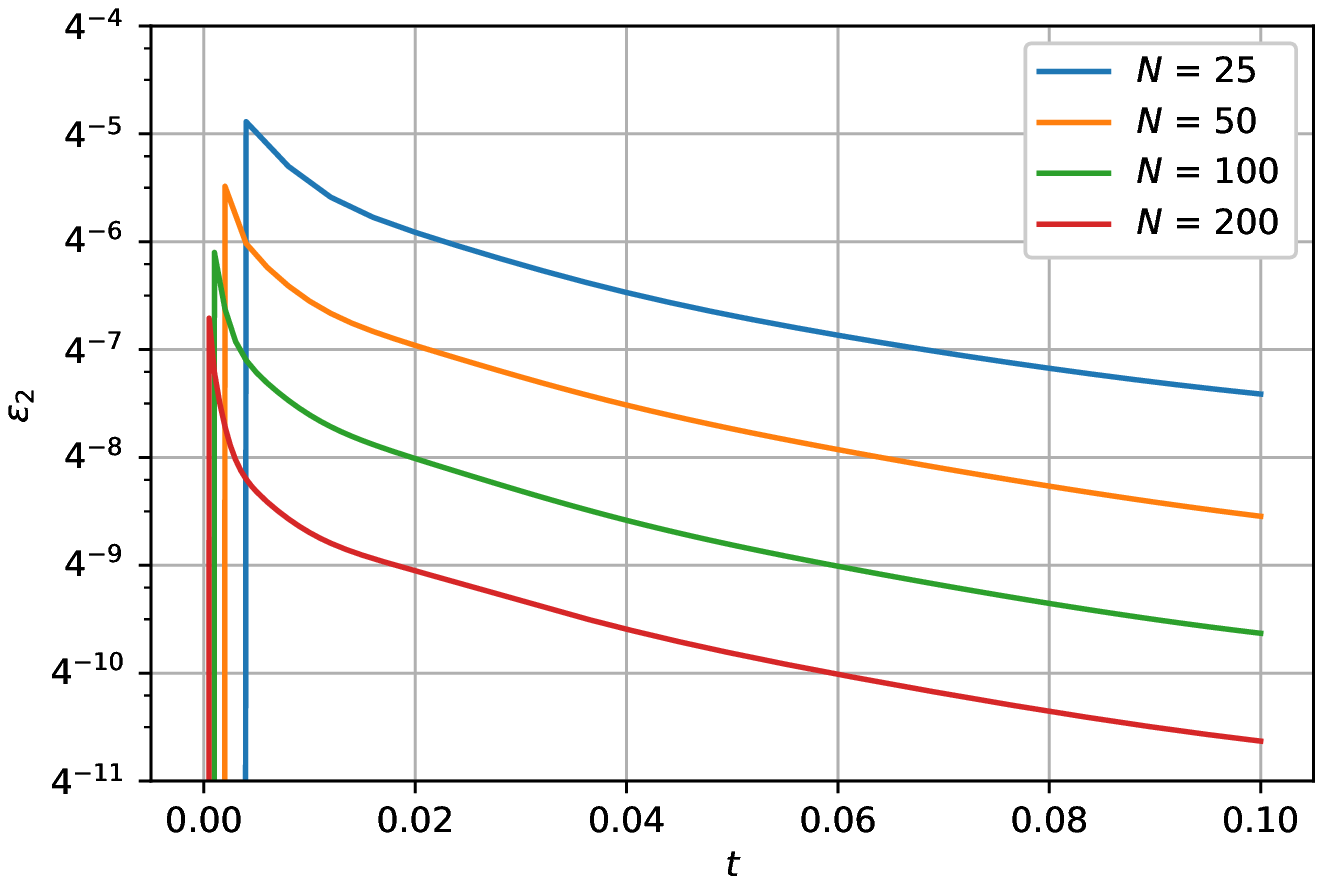} \includegraphics[width=0.45\linewidth]{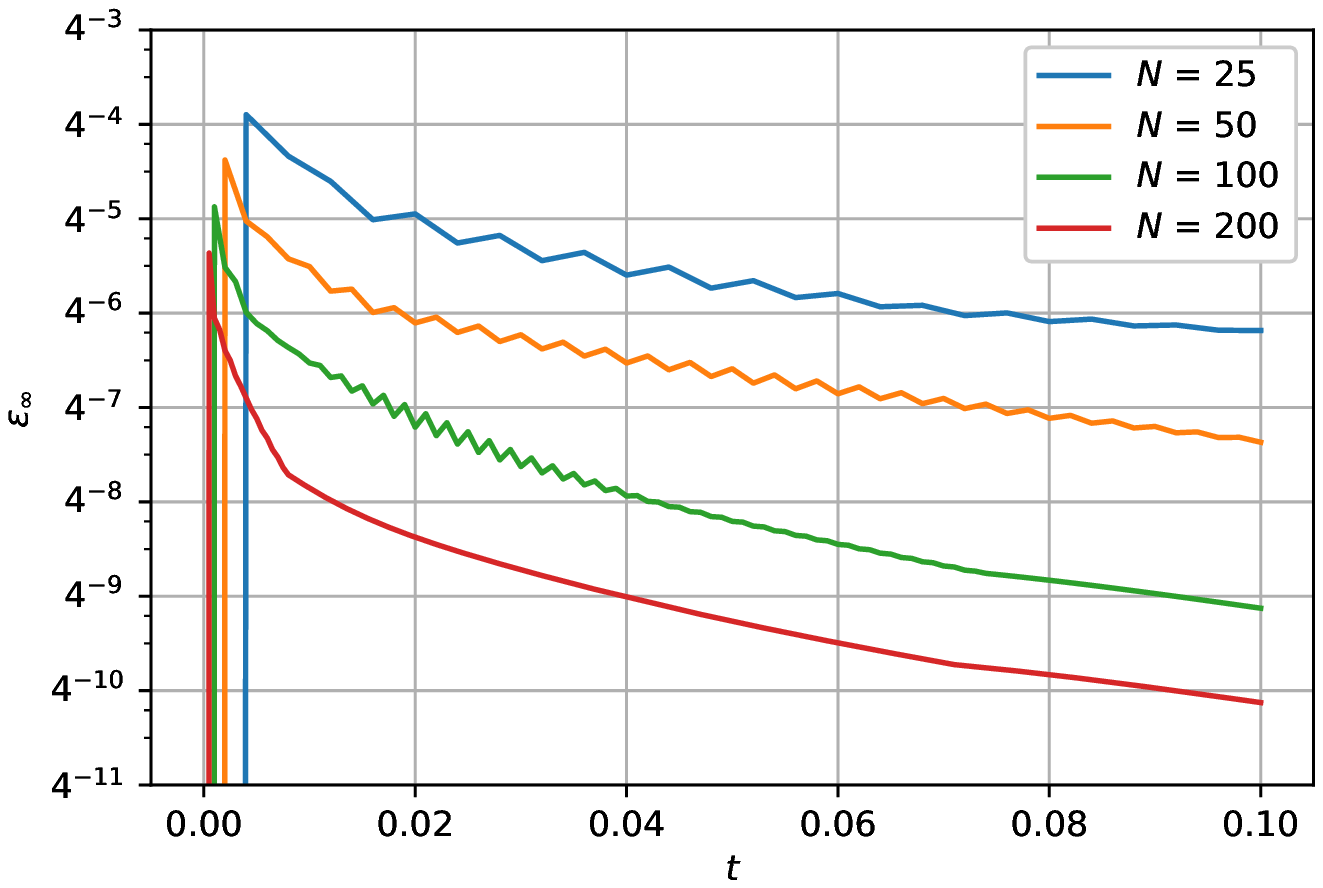} \\
\caption{The accuracy of the three-level scheme: $N_1 = N_2 = 64$ (top), $N_1 = N_2 = 128$ (bottom).}
\label{f-5}
\end{figure}

The solution of the test problem on the grid $N_1 = N_2 = 128$ for different moments of time is shown in Fig.\ref{f-3}.
The accuracy of the solution to the nonstationary problem (\ref{2.4}), (\ref{2.5}) is estimated by the absolute discrepancy at individual points in time:
\[
\varepsilon_2(t^n) = \|v (\bm x,t^n) - y^n(\bm x)\|,
\quad \varepsilon_\infty (t^n) = \|v (\bm x,t^n) - y^n(\bm x)\|_\infty ,
\quad n = 0, \ldots, N .
\]
Given that the accuracy of the decomposition scheme critically depends on the discretization over space, we perform calculations on two grids over space: $N_1 = N_2 = 64$ and $N_1 = N_2 = 128$.
The error of the approximate solution using different time steps $\tau = T/ N$ for the two-level scheme is presented in Fig.\ref{f-4}.
We observe a convergence of approximately the first-order in time. Note also that the accuracy is almost independent of the spatial grid.
Similar data for the three-level scheme are given in Fig.\ref{f-5}. We observe approximately second-order convergence on $\tau$.
When the mesh is thickened in space, the nonmonotonicity of the solution appears (see the data on $\varepsilon_\infty$).

\begin{figure}[htbp]
\centering
\includegraphics[width=0.5\linewidth]{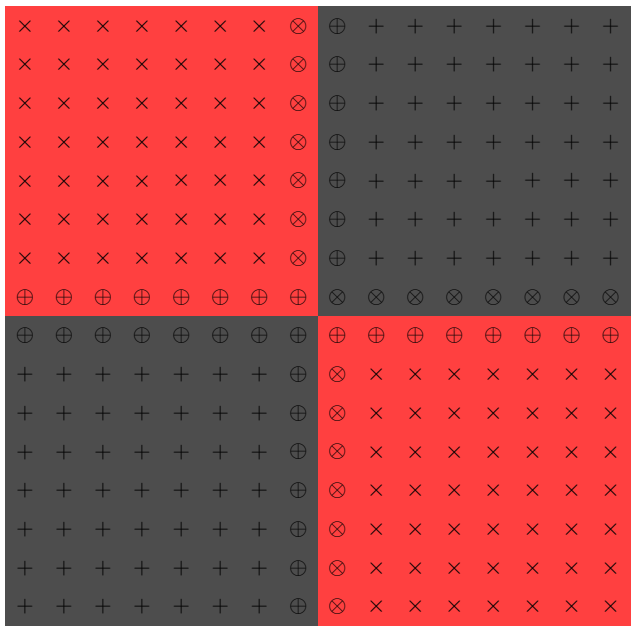} 
\caption{Red-black decomposition of the computational domain: $\bm  +$ --- nodes from $\omega_1$, $\bm \times$ --- nodes from $\omega_2$, $ \bigcirc$ --- exchange nodes.}
\label{f-6}
\end{figure}

We restrict ourselves to the case of decomposition into four non-overlapping subdomains (see Fig.\ref{f-6}).
The grid operator $A$ is set on a five-point template, and there is no exchange between the two pairs of subdomains.
This fact allows us to distinguish two enlarged subdomains; these subdomains correspond to red-black decomposition and the setting $p=2$.
When using decomposition schemes with the selection of the diagonal part $\bm A_0$ of the operator matrix $\bm A$, the transition to a new level by time is provided by independently solving four problems in subdomains in the original decomposition.

\begin{figure}[htbp]
\centering
\includegraphics[width=0.45\linewidth]{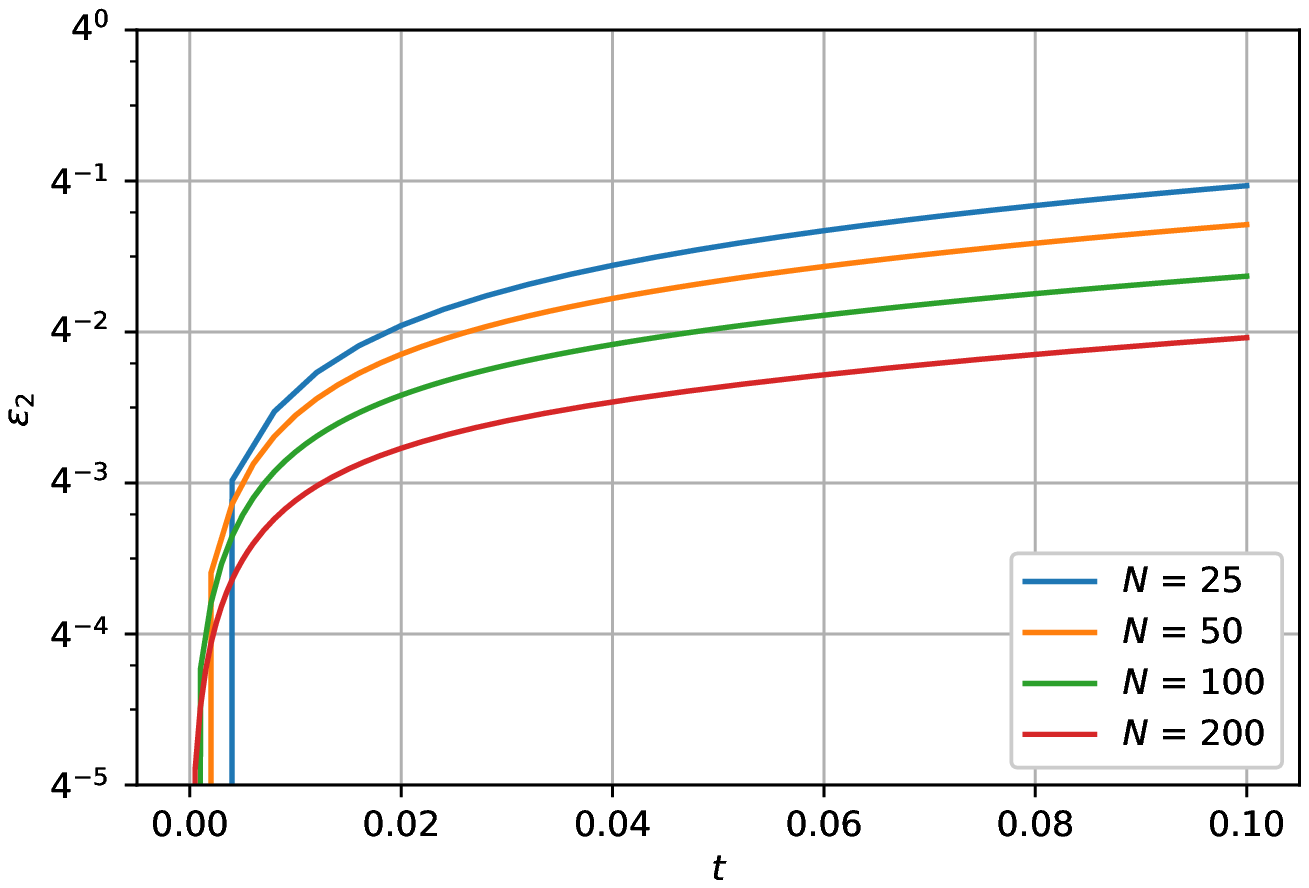} \includegraphics[width=0.45\linewidth]{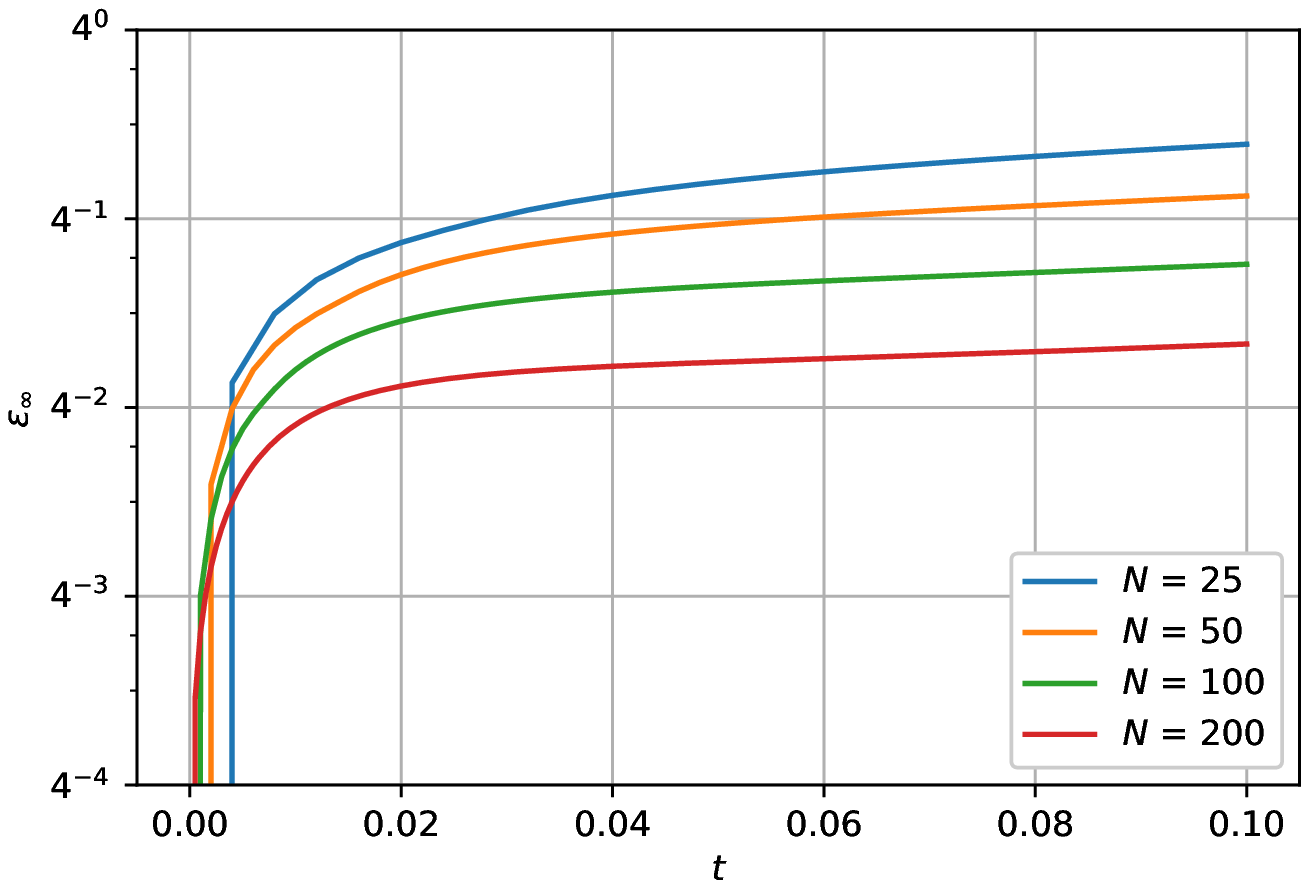} \\
\includegraphics[width=0.45\linewidth]{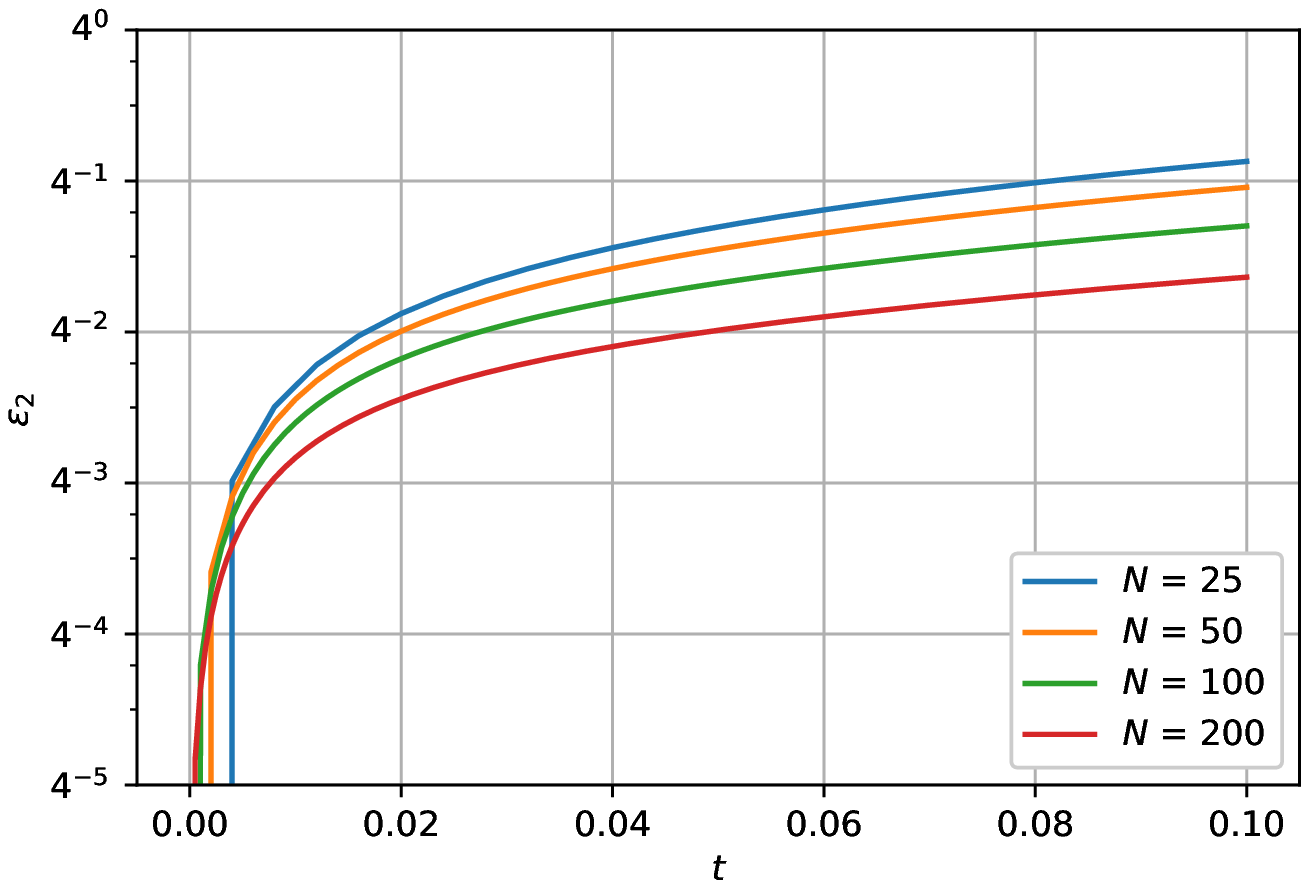} \includegraphics[width=0.45\linewidth]{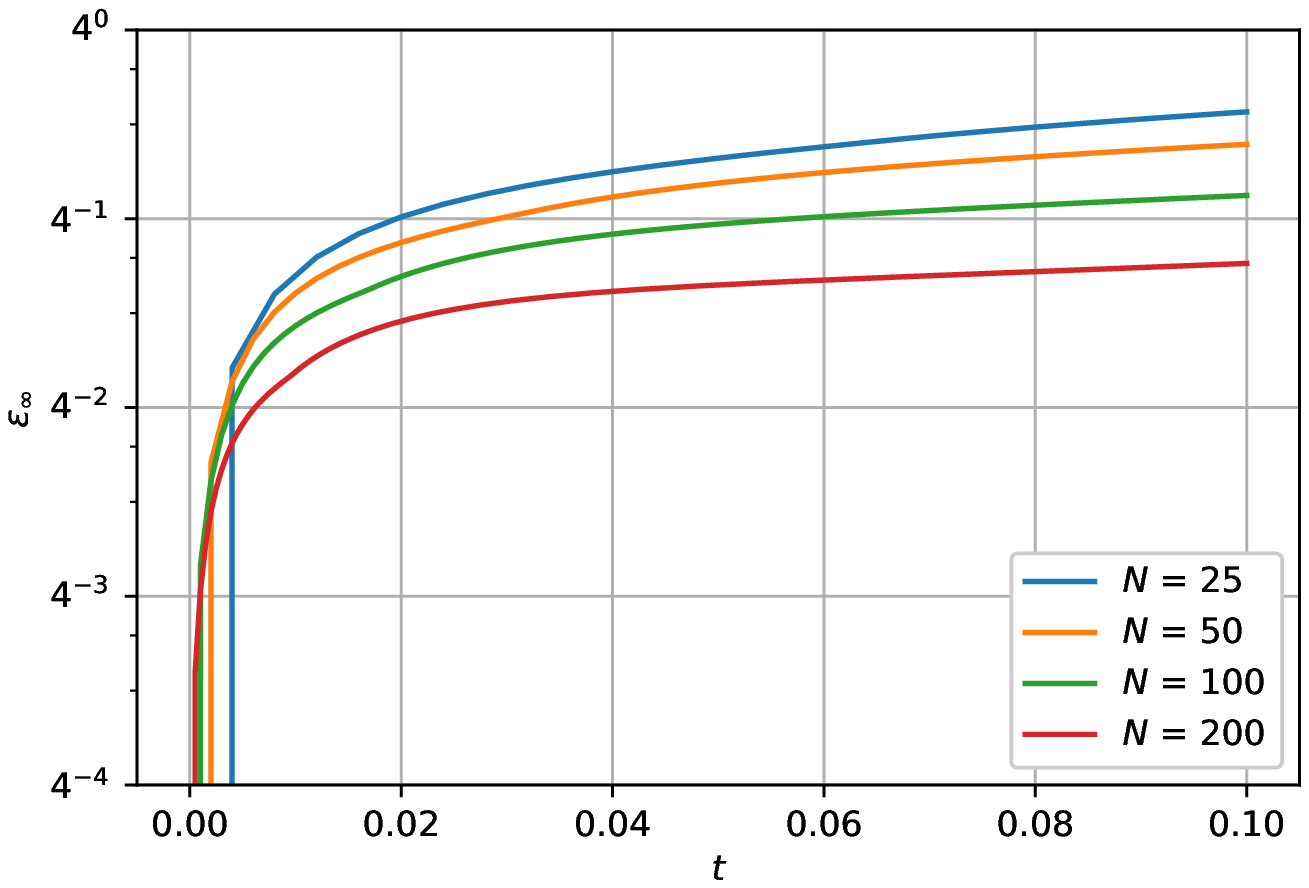} \\
\caption{The accuracy of the two-level domain decomposition scheme: $N_1 = N_2 = 64$ (top), $N_1 = N_2 = 128$ (bottom).}
\label{f-7}
\end{figure}

The accuracy of the two-level decomposition scheme (\ref{4.6}), (\ref{4.7}) with $\sigma = 1$ in the approximate solution of the test problem is shown in Fig.\ref{f-7}.
Comparison with the data in Fig.\ref{f-4} shows a significant drop in accuracy compared to the conventional implicit scheme (\ref{5.1}), (\ref{5.2}).
The domain decomposition schemes we are considering belong to the class of conditionally convergent schemes.
Because of this, when using a more detailed computational grid over the space, we find an approximate solution with less accuracy.

\begin{figure}[htbp]
\centering
\includegraphics[width=0.45\linewidth]{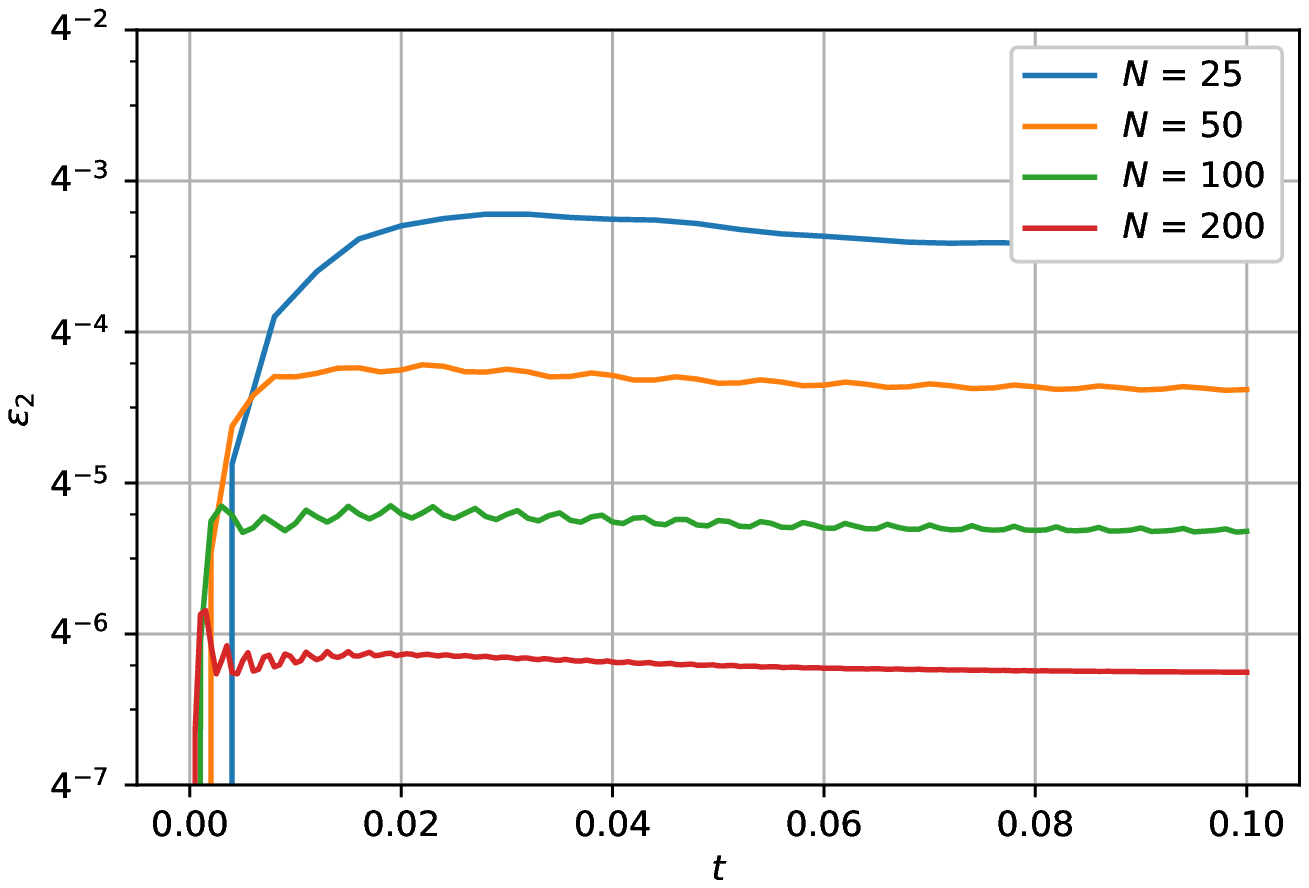} \includegraphics[width=0.45\linewidth]{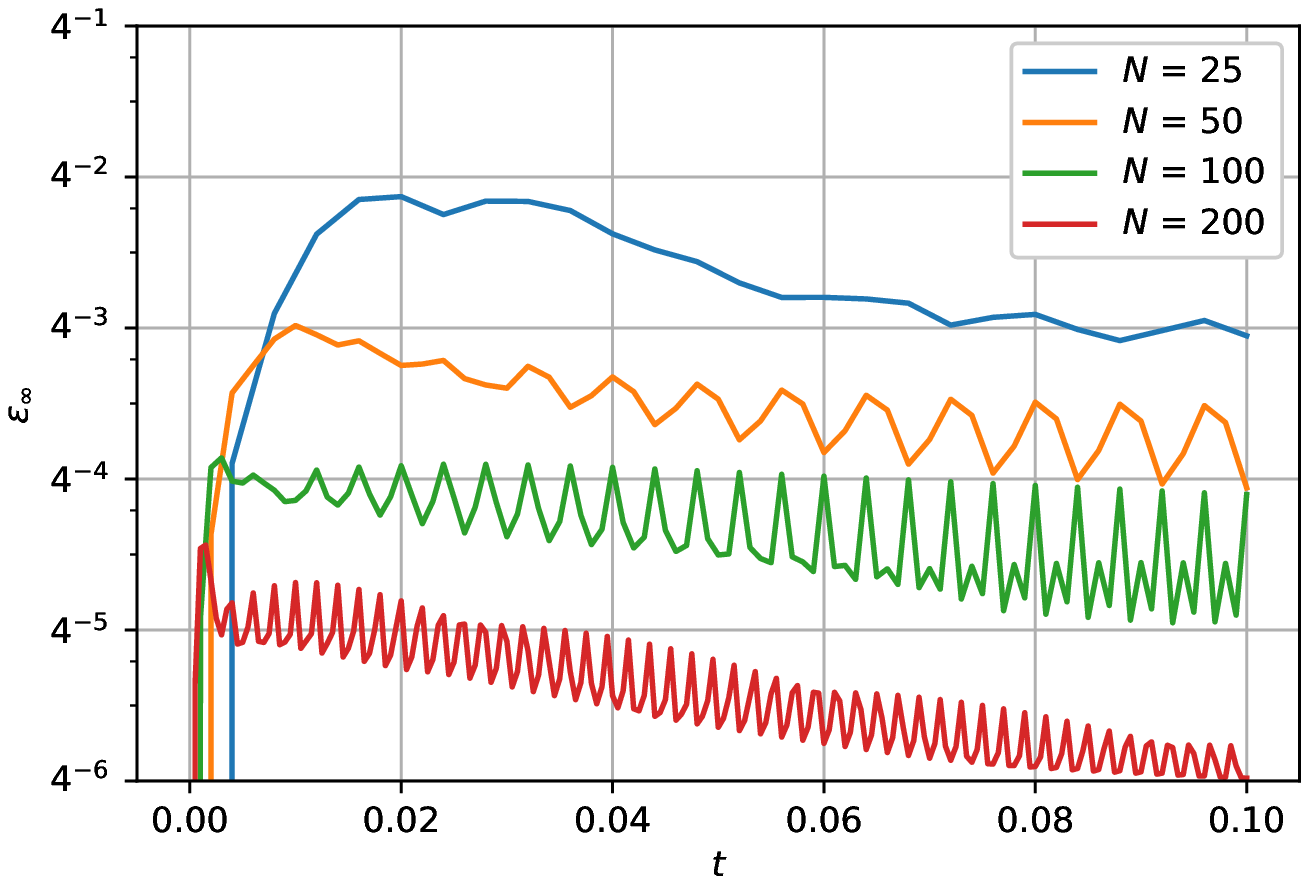} \\
\includegraphics[width=0.45\linewidth]{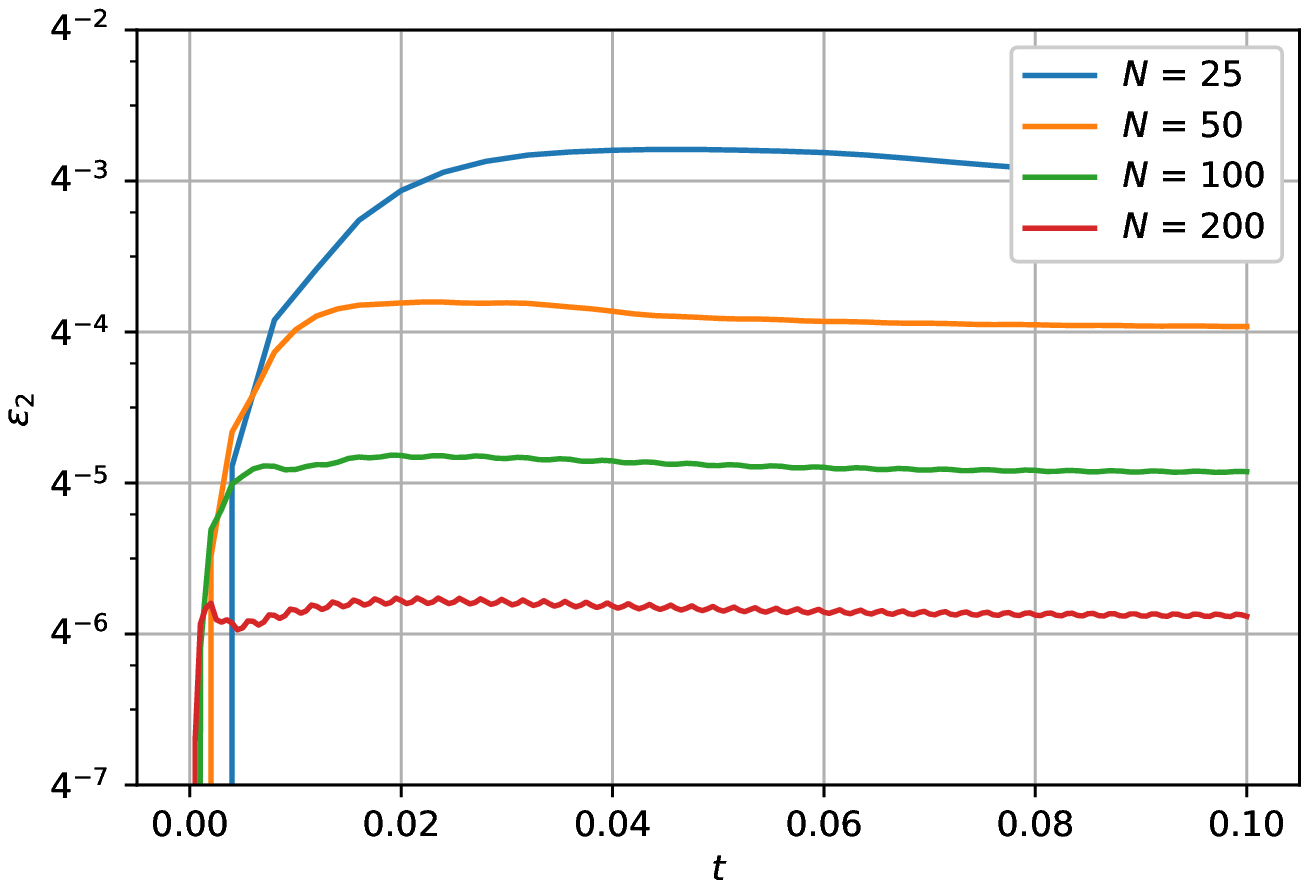} \includegraphics[width=0.45\linewidth]{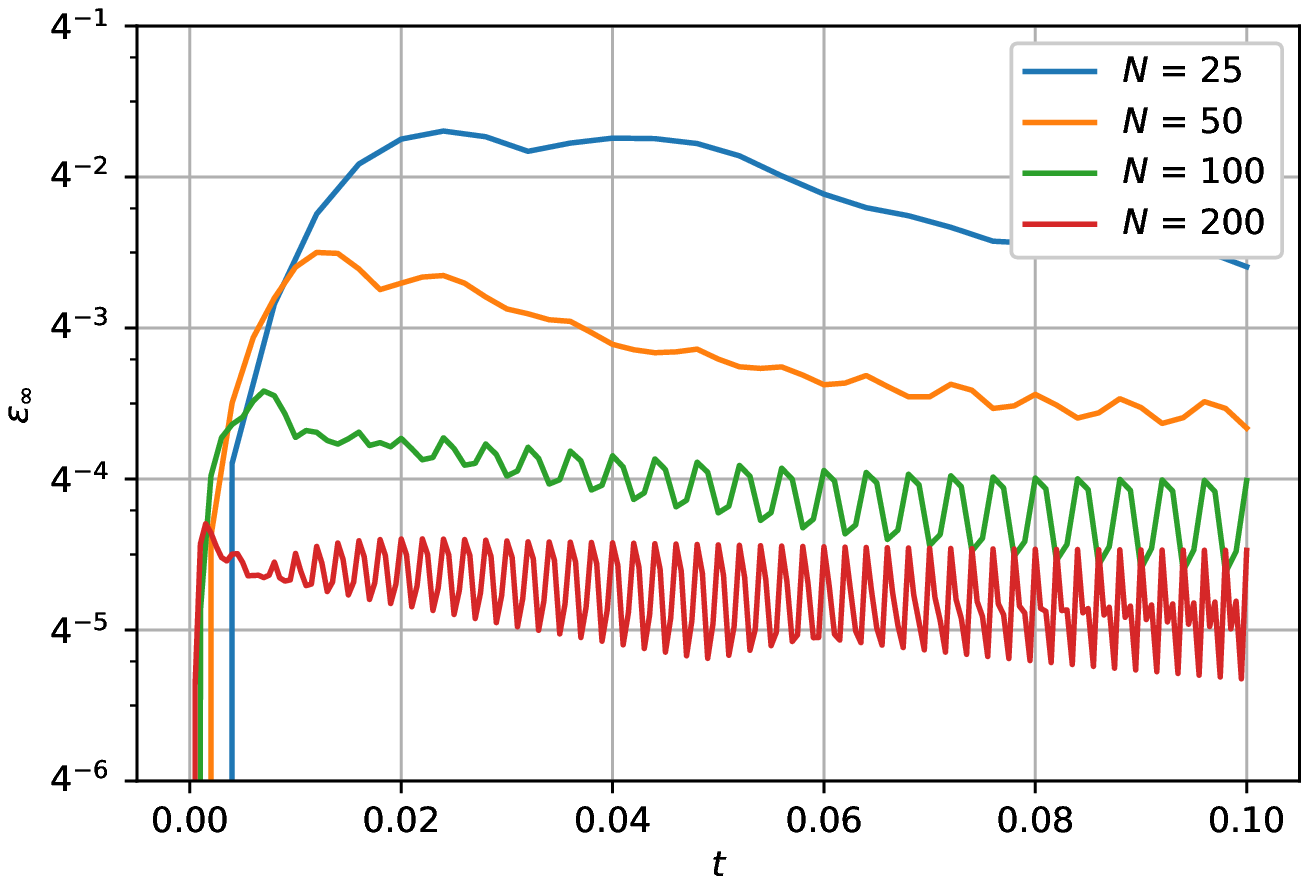} \\
\caption{The accuracy of the three-level domain decomposition scheme: $N_1 = N_2 = 64$ (top), $N_1 = N_2 = 128$ (bottom).}
\label{f-8}
\end{figure}

\begin{figure}[htbp]
\centering
\includegraphics[width=0.45\linewidth]{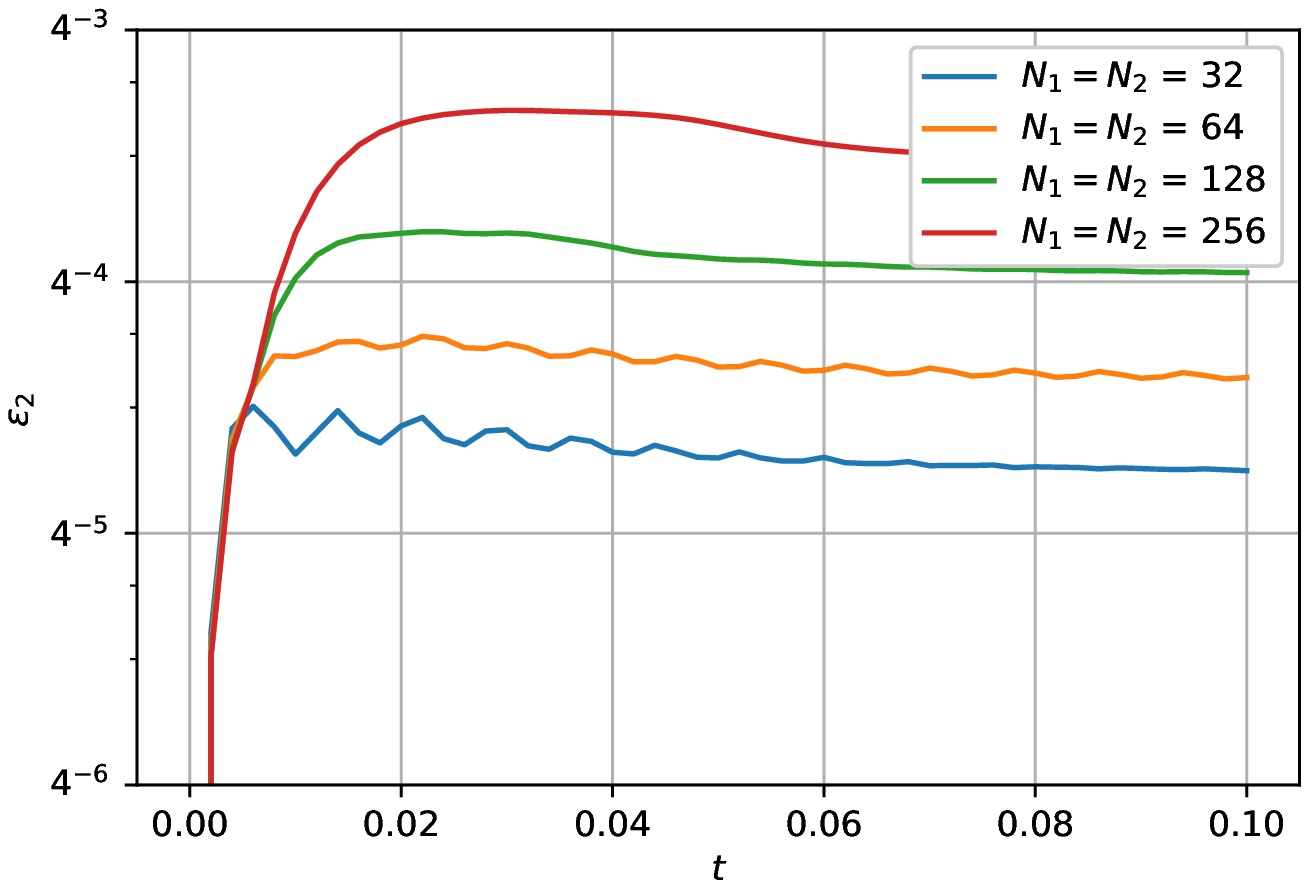} \includegraphics[width=0.45\linewidth]{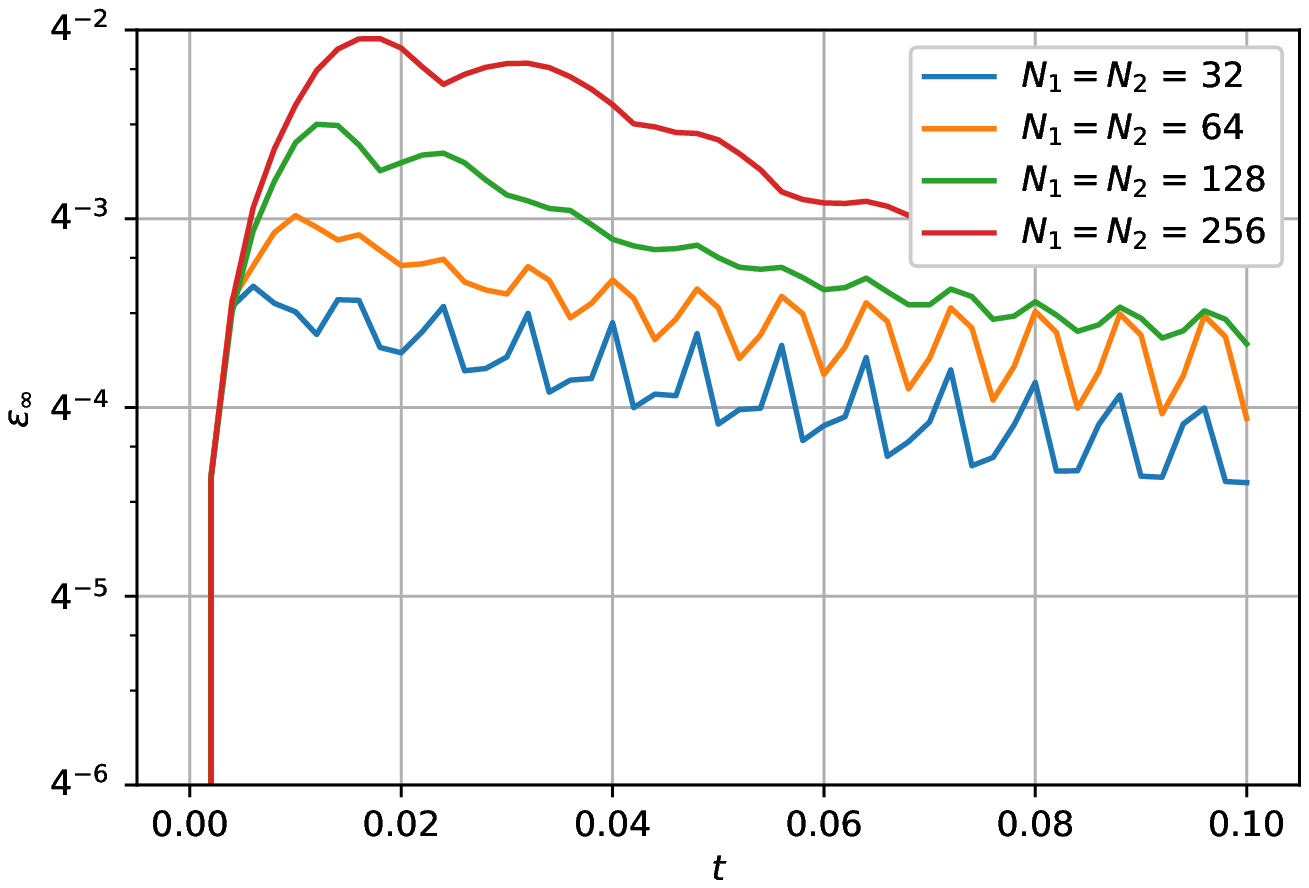} \\
\includegraphics[width=0.45\linewidth]{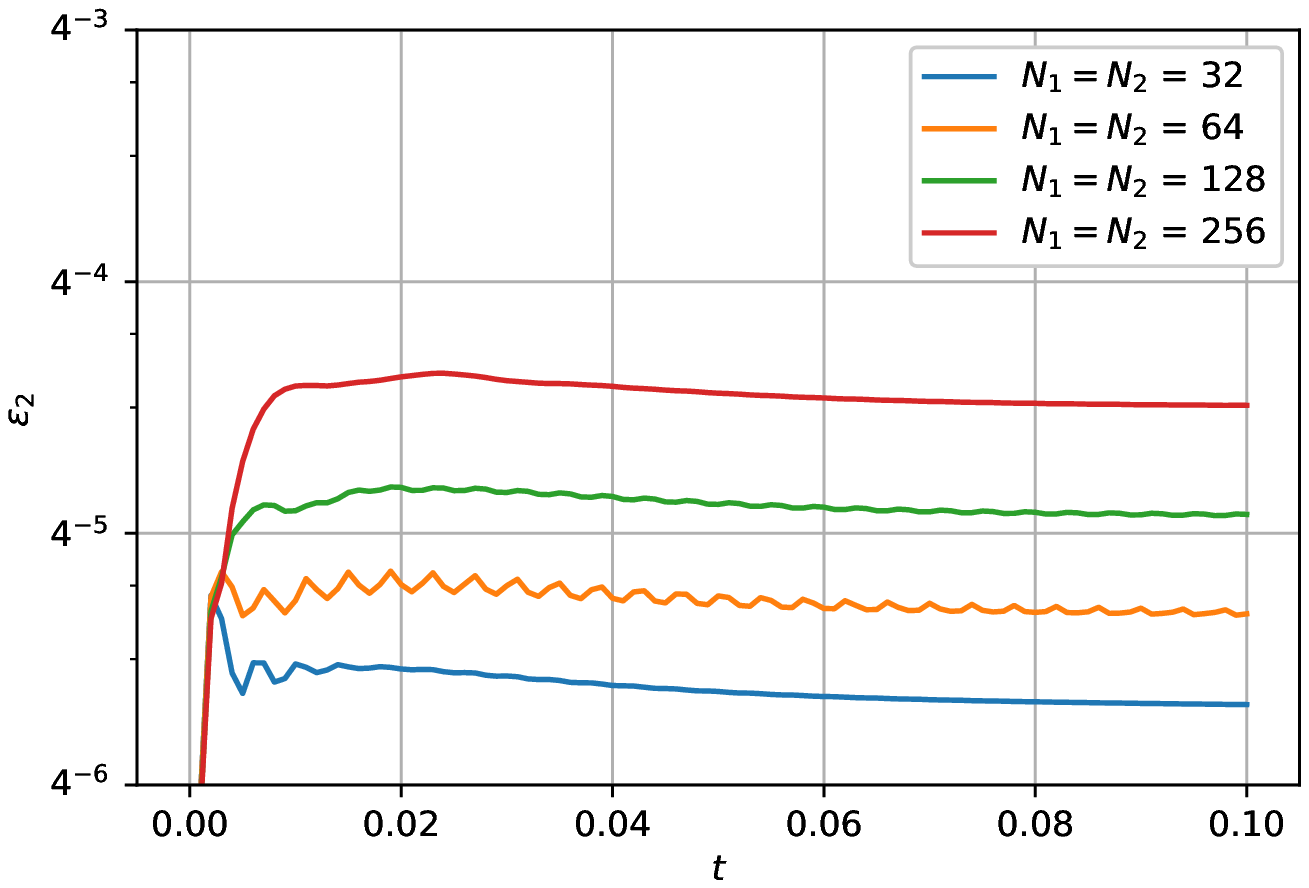} \includegraphics[width=0.45\linewidth]{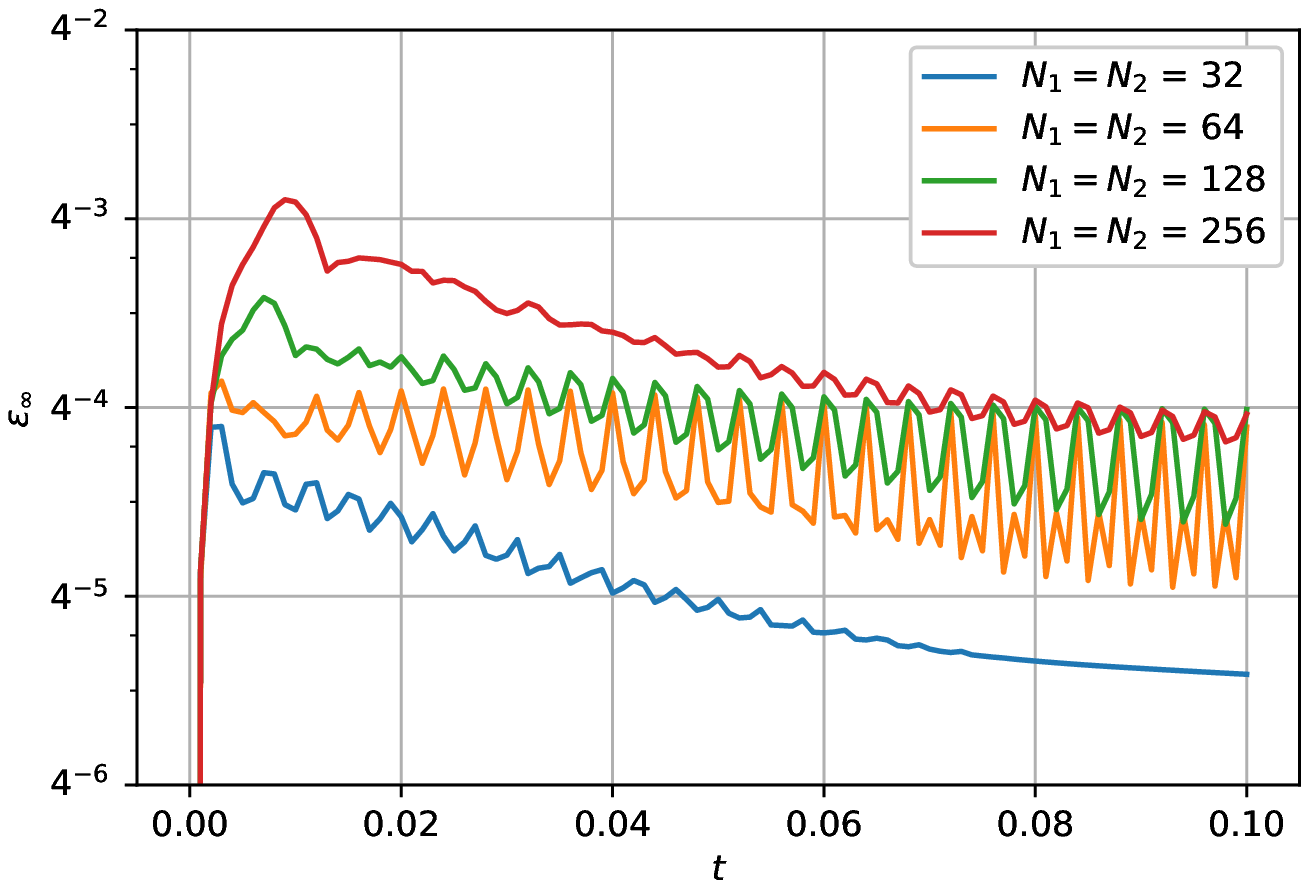} \\
\caption{The accuracy of the three-level domain decomposition scheme: $N = 50$ (top), $N = 100$ (bottom).}
\label{f-9}
\end{figure}

\begin{figure}[htbp]
\centering
\includegraphics[width=0.45\linewidth]{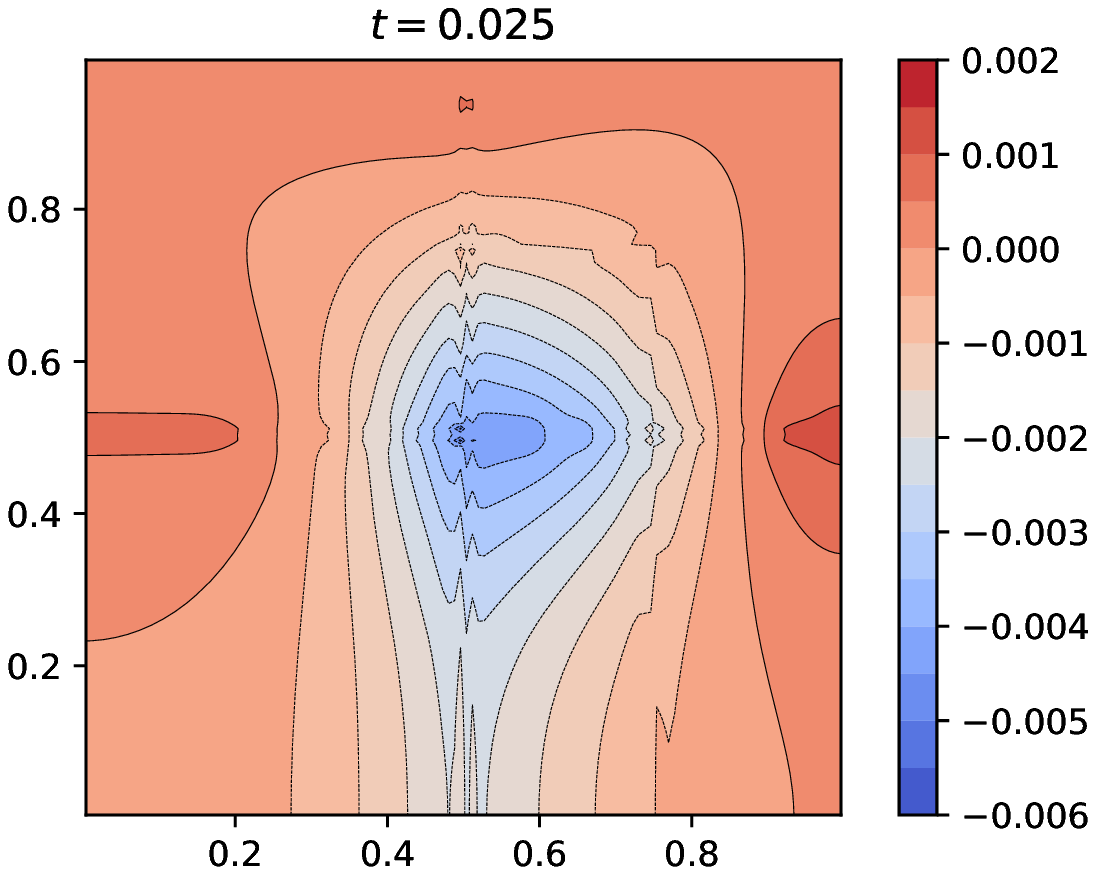} \includegraphics[width=0.45\linewidth]{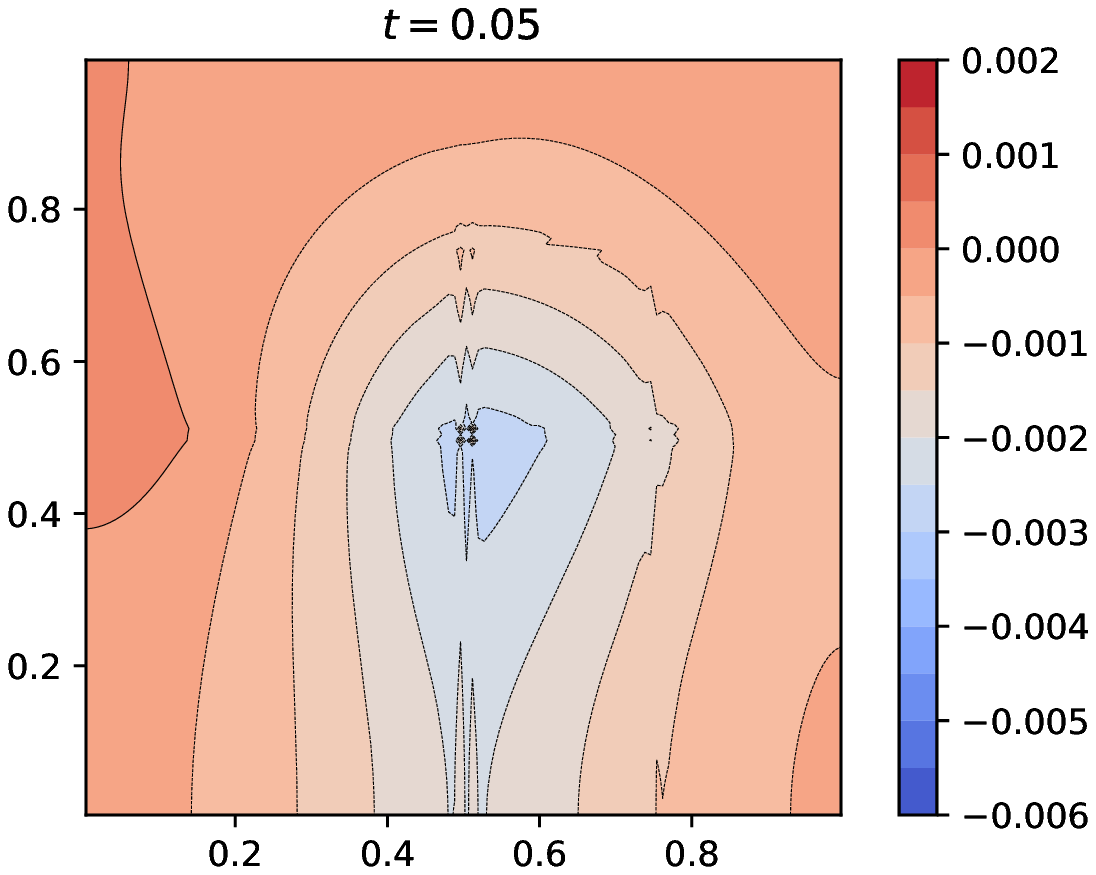} \\
\includegraphics[width=0.45\linewidth]{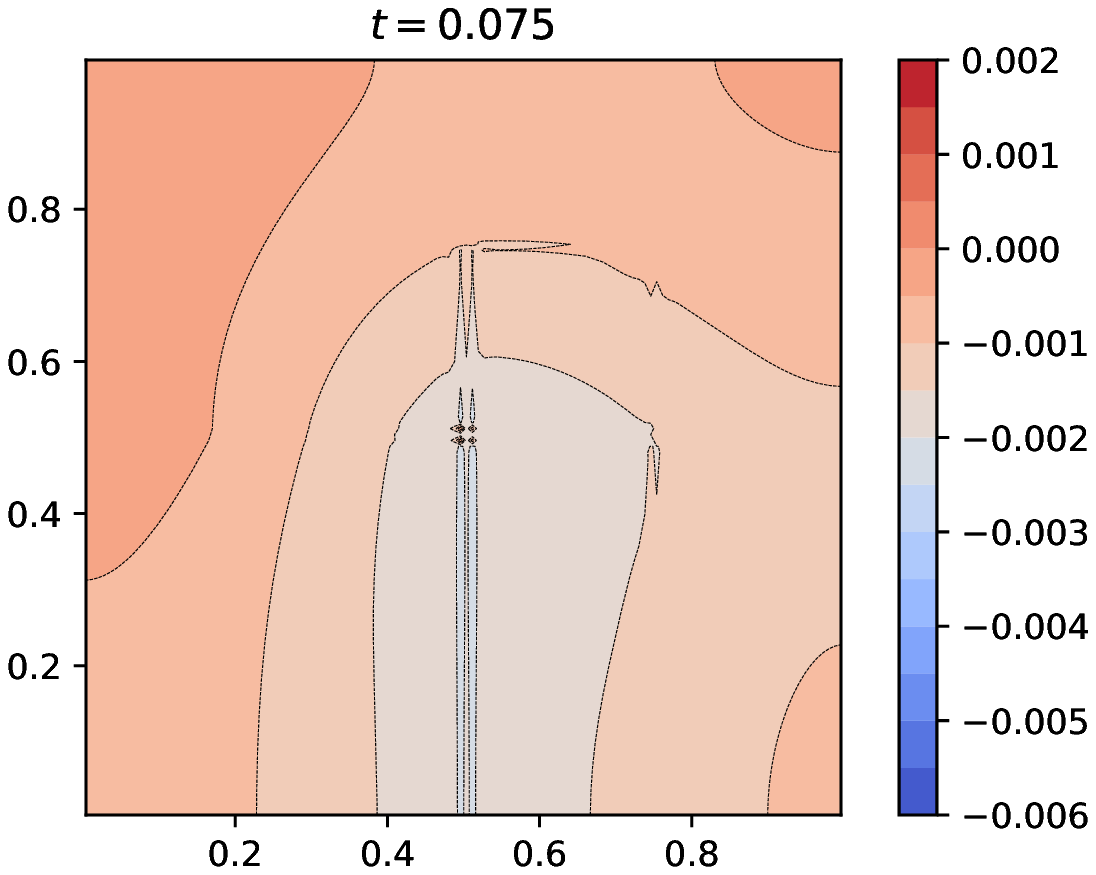} \includegraphics[width=0.45\linewidth]{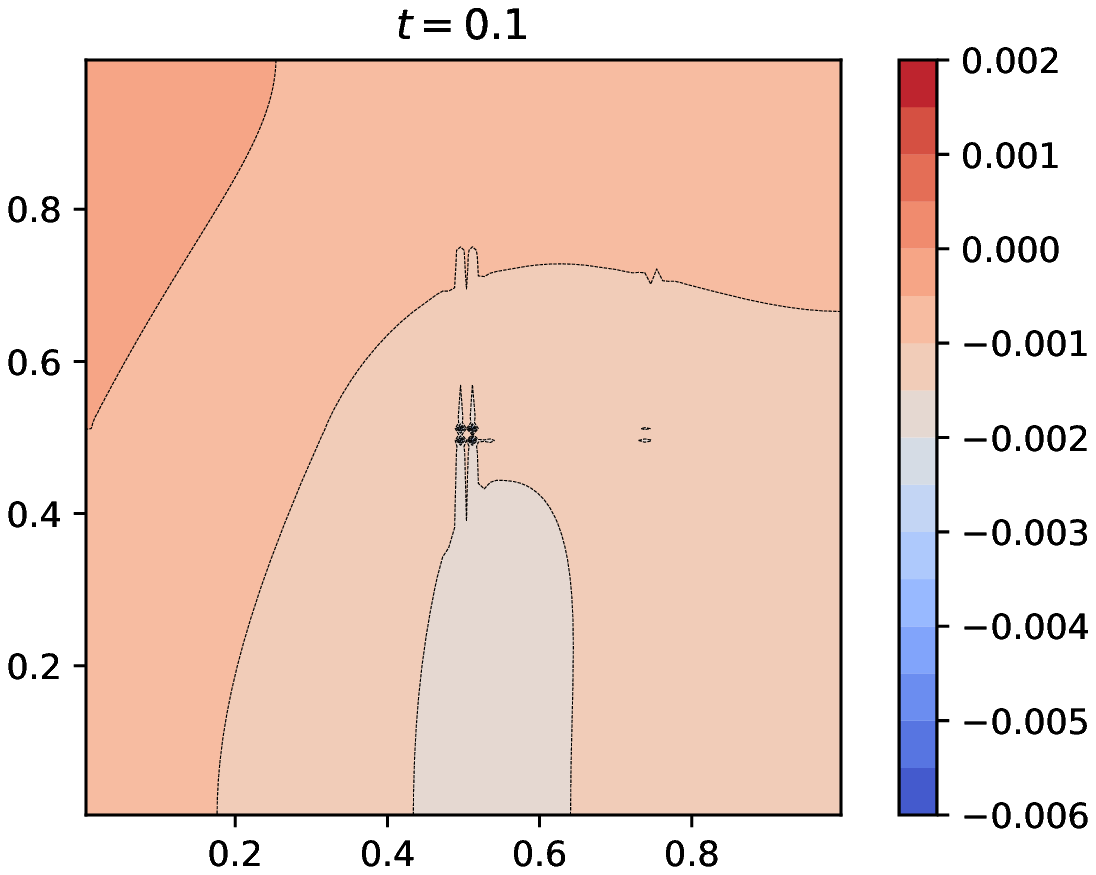} \\
\caption{The error of the three-level domain decomposition scheme: $N_1 = N_2  = 100$, $N = 100$.}
\label{f-10}
\end{figure}

The accuracy of the three-level decomposition scheme (\ref{4.15}), (\ref{4.16}) with $\sigma = 1/2$ is illustrated by Fig.\ref{f-8}.
Comparison with data from the standard three-level scheme (\ref{5.2})--(\ref{5.4}) in Fig.\ref{f-5} again demonstrates a decrease in the accuracy of the approximate solution when decomposing the domain.
The dependence of the error due to domain decomposition is most powerfully demonstrated when using different meshes by space (Fig.\ref{f-9}).
The error of the approximate solution $v(\bm x, t^n) - y^n(\bm x), \ \bm x \in \omega$ at some time $t = t^n$ is shown in Fig.\ref{f-10}.
The errors are localized near the boundaries of the subdomains.
The main conclusion is that the three-level domain decomposition scheme gives an approximate solution with significantly higher accuracy than the two-level domain decomposition scheme.
This statement is entirely consistent with the results of our theoretical study.

\section{Conclusions}\label{sec:6}

\begin{enumerate}
\item A typical boundary value problem for a parabolic equation in a rectangle is considered.
After the finite-difference approximation in space, the Cauchy problem for the first-order differential-operator equation is formulated.
A standard estimation of the stability of the approximate solution by initial data and the right-hand part for unconditionally stable two-level schemes with weights is given. It is necessary to simplify the grid problem on a new time level.
\item The decomposition of the solution to the Cauchy problem based on a particular choice of restriction and extension operators is performed.
In this case, we obtain a coupled system of equations for the solution components.
An explicit representation of restriction and extension operators in the decomposition of the grid domain with and without overlapping subdomains is given.
\item We obtain stability conditions for two- and three-level schemes in the corresponding Hilbert spaces in the solution decomposition when: (i) the diagonal part of the operator matrix of the problem is extracted, and (ii) the triangular splitting of the operator matrix is used.
An analysis of the accuracy of the proposed iterative-free domain decomposition schemes is performed.
\item Theoretical results are illustrated by data of numerical experiments for the test problem.
In particular, the best prospects for the practical use of three-level schemes of the second order of approximation by a time when the diagonal part of the operator matrix of the problem is extracted on a new time level are confirmed.
\end{enumerate} 


\end{document}